\documentstyle[amssymb,11pt]{article} 
\newcommand{\RR}{\mathbb R}  
\newcommand{\NN}{\mathbb N} 
\newcommand{\CC}{\mathbb C} 
\newcommand{\ZZ}{\mathbb Z} 
\newcommand{\del}{\partial}  
\newcommand{\e}{\varepsilon}  
  
\newcommand{\ha}{\frac{1}{2}} 
  
\newcommand{\Lap}{\Delta}  
  
\newcommand{\hZ}{\widehat{Z}}  
\newcommand{\cin}{{\cal C}^\infty} 
\newcommand{\Vb}{{\cal V}_b} 
\newcommand{\ran}{\mbox{\rm ran\,}} 
\newcommand{\ind}{\mbox{\rm ind\,}} 
\newcommand{\Ind}{\mbox{\rm Ind\,}} 
\newcommand{\nul}{\mbox{\rm null\,}} 
\newcommand{\spec}{\mbox{\rm spec}} 
\newcommand{\Id}{\mbox{\rm Id}}
\newcommand{\id}{\mathrm Id} 
\newcommand{\ff}{\mathrm ff} 
\newcommand{\lb}{\mathrm lb}  
\newcommand{\rb}{\mathrm rb}  
\newcommand{\ls}{\mathrm ls} 
\newcommand{\rs}{\mathrm rs} 
\newcommand{\ds}{\mathrm ds} 
\newcommand{\db}{\mathrm db} 
\newcommand{\hs}{\mathrm hs} 
\newcommand{\Diff}{\mbox{\rm Diff\,}} 
\newcommand{\Spi}{S\kern-6pt/} 
\newcommand{\Dir}{D\kern-6pt/} 
\newcommand{\DDir}{{\cal D}}
\newcommand{\Ima}{\mbox{\rm Im}} 
\newcommand{\Lam}{\Lambda^{\mathrm max}} 
\newcommand{\Ker}{\mbox{\rm Ker\,}} 
\newcommand{\coker}{\mbox{\rm coker\,}} 
\newcommand{\Hom}{\mbox{\rm Hom\,}} 
\newcommand{\ilg}{\mbox{\rm ilg\,}}
\newcommand{\Ls}{\mathrm Ls}

\newtheorem{theorem}{Theorem}

\begin{document}  
  
\title{Dirac operators, heat kernels and microlocal analysis \\  
Part II: Analytic surgery}   
\author{Rafe Mazzeo\thanks{Supported by the NSF with a Young Investigator   
Fellowship and by grant \#DMS-9626382}\\ Stanford University \\   
mazzeo\@@math.stanford.edu \and Paolo Piazza \\ Universit\`a di Roma  
``La Sapienza'' \\ piazza\@@mat.uniroma1.it}    

\maketitle  

\noindent
{\bf Abstract.} Let $X$ be a closed Riemannian manifold and let
$H\hookrightarrow X$ be an embedded hypersurface. Let $X=X_+ \cup_H X_-$
be a decomposition of $X$ into two manifolds with boundary, 
with $X_+ \cap X_- = H$.
In this expository article, surgery -- or gluing -- formul\ae\ 
for several geometric and spectral invariants associated to 
a Dirac-type operator $\eth_X$ on $X$ are presented.
Considered in detail are: the index of $\eth_X$, the index bundle and the 
determinant bundle associated to a family of such operators, the eta invariant
and the analytic torsion.  In each case the precise form of the surgery
theorems, as well as the different techniques used  
to prove them, are surveyed.\\

  
\section{Introduction}  

The behaviour of global invariants for Dirac operators and Laplacians 
with respect to decompositions of their underlying compact Riemannian 
manifolds has become a topic of much interest over the past several
years. We are thinking here of {\it geometric} invariants such
as index and determinant bundles, as well as {\it spectral} invariants
such as the eta invariant and analytic torsion
So-called `gluing theorems' for these invariants provide 
new insights into their nature and have facilitated their use in other 
areas. One technique to study these problems was developed by the
first author and Melrose \cite{asatei}, and McDonald \cite{Mc}, and 
is called {\it analytic surgery}. In this paper we give a brief 
introduction to this method and to a few of the problems for which it 
has proved useful, and also to survey a few other methods developed 
by other authors to study gluing problems. 

In the most general terms, suppose that we are given a decomposition
of the compact manifold $X$ into two pieces, $X = X_+\cup X_-$, where
$X_\pm$ are submanifolds with boundary. Any geometric datum, such as a 
Riemannian metric $g$, a bundle $E$, or a spin structure and its associated 
Dirac operator $\eth$, restricts to give the corresponding structures on 
each of these pieces. In the next section we shall give precise definitions of 
some of the global invariants in which we are interested; 
for the sake of being concrete, and referring to that 
section for its definition, let us consider the eta invariant of a Dirac 
operator, $\eta(\eth)$. Setting aside, for the moment, the issue of boundary 
conditions, the simplest formulation of one of the problems we wish to 
discuss is whether  
there is a reasonable formula for $\eta(\eth_X)$ in terms of 
$\eta(\eth_{X_\pm})$; here, for any manifold $Z$, possibly with boundary, 
$\eth_Z$ denotes its Dirac operator relative to some fixed spin
structure and metric. Amongst the various considerations we shall need to 
address, even just to formulate a reasonable conjecture more precisely,
is the issue of boundary conditions, and also such matters as the dependence
of the eta invariant on the underlying metric. Upon doing this, it
will become apparent that it is very natural, or at least very
convenient, to study families of degenerating metrics, or families
of boundary conditions, and that the defect between $\eta(\eth_X)$
and $\eta(\eth_{X_\pm})$ is reasonably gauged by some measure
of the variation in these families. The problem then will be to
express this defect in some explicit way.

The structure of this paper will be somewhat informal, inasmuch
as we consider these various problems in successively greater
degrees of precision. In the rest of this introduction we formulate
the surgery problem somewhat more carefully; the reader should
note that different authors describe it in seemingly quite different 
ways, depending on their precise contexts and the applications 
they have in mind. We are trying here to present these approaches
from a more uniform perspective. 
After this `first-level' formulation, we proceed
to discuss two competing points of view related to the issue of 
boundary conditions: namely, is it more natural to consider
various geometric structures on a manifold with boundary $Z$
as smooth up to the boundary (and possibly having some product 
structure near the boundary $Y = \del Z$), or else as defined
relative to an infinite cylindrical geometry near the ends.
As was recognized already in \cite{APS}, these two points of view 
are essentially equivalent, but choosing one or the other
as primary tends to inform our intuitions in different ways.
(Of course, there are many other geometries on a manifold
with boundary, but these two have, up until now at least, played 
the most prominent r\^oles in the sorts of problems we consider.)  
We include here a brief overview of the calculus of $b$-pseudodifferential 
operators, as a preamble to the surgery calculus discussed later. 
We conclude the introduction by discussing the three principle methods 
used to study gluing theorems: those developed by Bunke, Vishik, and 
that of the first author, Melrose and Hassell.

In the remainder of the paper we shall, as promised, give more
careful explanations of many of these issues. In \S 2, we discuss
some of the different settings and invariants for which each of these
methods has proved useful, or at least has been applied. 
In the three succeeding sections we give more careful discussions of these 
three methods, concentrating, it must be admitted, on the final one. 
Unfortunately, we do not have the time or space to go too deeply into any of 
the analytic subtleties in any of these approaches, but instead wish to 
present them side by side, indicating some of their relative strengths and 
weaknesses in hopes that this will be useful for future applications.
In the final sections we give a more extended discussion of two
further applications of the surgery calculus: the first is to the
signature formula on manifolds with corners as in \cite{stmwc}, while 
the second is to gluing formul\ae\ for determinant bundles as in \cite{P3}.
 
\subsection{The surgery problem} 

We start, as above, with the decomposition 
\begin{equation}
X = X_+ \cup X_-, \qquad \mbox{where} \qquad X_+ \cap X_- = H
\label{eq:1.1}
\end{equation}
is a smooth, oriented hypersurface in $X$ and the pieces $X_\pm$ 
are smooth manifolds with boundary. (We are implicitly assuming that 
$H$ disconnects $X$; this is not necessary, but we consider only this 
case so as to minimize notation.) A Riemannian metric $g$ on $X$ induces 
metrics $g_\pm$ on $X_\pm$, and if $\Delta_Z$ denotes the (scalar) 
Laplacians on any one of these manifolds, $Z = X, X_\pm$, then a primitive 
form of the analytic surgery problem is to determine the relationship 
between ${\rm spec\,}(\Delta_X)$ and ${\rm spec\,} (\Delta_{X_\pm})$. 
(To define the latter quantities, we use, for example, Dirichlet 
conditions on $H$.) While precise relationships between individual eigenvalues 
are generally impossible to establish, it is easier to find relationships 
between some aggregate invariants of these spectra, such as the 
determinants $\det(\Delta)$, or even between their resolvents
$(\Delta - \lambda)^{-1}$ or heat kernels, $\exp(-t\Delta)$.

At a slightly higher level of complexity, suppose that $\eth_X$ is the
Dirac operator on $X$ with respect to some fixed spin structure and the 
metric $g$, or even simply a generalized Dirac-type operator (which does 
not require a spin structure {\it per se}). The issue of boundary conditions 
for the restricted Dirac operators $\eth_{X_\pm}$ on $X_\pm$ is now more 
subtle, and it is well-known that one must use global boundary conditions, 
of the sort introduced by Atiyah, Patodi and Singer, to obtain an elliptic 
boundary problem.  We discuss this in the next subsection. At any rate, 
having done this, once again we ask for relationships between the spectra
of these operators or between the resolvents or heat kernels associated
to  their squares, $\eth^2$. The most common global spectral invariants
in this context are the eta invariant $\eta(\eth)$, which is of particular 
interest only when $\dim X$ is odd, and the analytic torsion,
$\tau(X)$.
Again referring to the eta invariant for concreteness, it is not quite true
in general that the eta invariant for $\eth_X$ is simply the sum of the 
eta invariants for $\eth_{X_\pm}$. So, as indicated earlier, the problem 
reduces to finding a good formula for the 
defect between these two expressions (just as the eta invariant is the 
defect between the two sides of the index formula).  From this
point of view, the defect arises because the boundary conditions
which the $\eth_{X_\pm}$ inherit by restriction from $\eth_X$
do not match the natural global APS boundary conditions on $X_\pm$. 
One may change perspective, though, and consider instead the
eta invariant of $\eth_X$ relative to a family of metrics $g_\e$ 
on $X$ which degenerate along the hypersurface $H$. Denoting this
family of operators by $\eth_{X,\e}$, then $\eta(\eth_{X,\e})$ depends on this family 
of metrics, but in an understandable way, at least for $\e > 0$ (here $\e =0$
corresponds to some sort of degenerate limit). Now the problem
becomes to determine the defect between the eta invariant
of the limiting operator and the limit of the eta invariants. 
One may pose a similar problem for the log of the analytic
torsion, $\log \tau(X)$. 

The types of metric degenerations we shall discuss, and shall
call {\it surgery degenerations}, arise when $g_\e$ elongates 
transversally to $H$, but stays bounded (and converges smoothly 
to a limit) away from $H$ in such a way that in the limit as $\e 
\rightarrow 0$, the interiors of the components $X_\pm$ inherit 
complete metrics with asymptotically cylindrical end structures. 
It may indeed be reasonable, and even better in some
circumstances, to study other sorts of degenerations.
For example, there is an enormous literature concerning
degenerations of compact Riemann surfaces endowed with
their hyperbolic metrics into limits with hyperbolic cusp ends --- 
indeed, a dense set of points on the boundary of
Teichm\"uller space of a surface may be `reached' in this way ---
and Arakelov degenerations have also received considerable 
attention, cf. \cite{Wen}, \cite{Jost}. 
Whether these other geometries are more
favourable for some of our spectral questions is not known.

Beyond the questions concerning these numerical invariants
are some others, particularly when one studies families of 
degenerating metrics and their associated Laplacians or Dirac 
operators. For the cylindrical degenerations we shall study, the 
spectrum of $\Delta_{g_\e}$ or $\eth_{g_\e}$ is discrete when 
$\e > 0$, but is continuous when $\e = 0$ (with possibly
some additional discrete spectrum). 
It is natural to enquire how the transition 
between these two states takes place; in particular, how accurately
may one describe the convergence of the discrete spectrum to continuous 
spectrum.  This sort of question seems relevant in light of the 
extensive recent work on Novikov-Shubin invariants, which are some 
sort of measure of the `germ' of continuous spectrum at $0$ for the 
Hodge Laplacian on differential forms on universal covers of compact
manifolds.  This sort of geometry is quite different
from cylindrical end geometry, but it is clear that there is much to 
be learned about the fine structure of the continuous spectrum of 
geometric operators on complete manifolds. 

\subsection{\mbox{\rm APS} vs. $L^2$ boundary conditions} 
Obviously there is a substantial geometric difference between complete 
metrics (with cylindrical ends) on $X_\pm$ and the restrictions of the 
metric $g$ to these components, and just as obviously there are substantial 
analytic differences between the Laplacians or Dirac operators for these 
metrics.  Whether to approach the surgery problem via restriction or 
degeneration of metrics is a moot point: as we shall see, each point 
of view has its strengths and weaknesses. But at heart is the purely 
qualitative and subjective question concerning which classes of functions 
or metrics or differential operators on a compact manifold with boundary 
one should consider to be the most natural.

At first, it seems odd to say anything other than that the classes of 
objects, e.g. functions, metrics, etc., which are smooth up to the 
boundary in the usual sense are the most natural ones. However, there 
is a good argument to be made that this is not necessarily the case. 
At the very least, the class of metrics with asymptotically cylindrical 
ends, the geometric elliptic operators corresponding to these metrics,
and finally, the class of $b$-pseudodifferential operators which 
generalize these, have proven to be essential in a 
number of recent geometric investigations. Amongst these we wish to 
mention the recent `direct' proof of the index theorem of Atiyah, 
Patodi and Singer on manifolds with boundary obtained 
by Melrose \cite{tapsit}. Although in many ways equivalent 
to the original approach to this theorem, its slightly different
perspective and the use of $b$-pseudodifferential operators led the 
way to previously unknown results, such as the proper generalization of the
index theorem for families of Dirac operators on manifolds with boundary 
by the second author and Melrose \cite{MP2}, \cite{MP3}, following partial 
results by Bismut and Cheeger \cite{BC} using an older approach. 
The predecessor to this paper \cite{P1} surveys these results.  
In addition, the long exact 
sequence in analytic $K$-homology for manifolds with boundary,
or even with corners, had proved 
somewhat difficult to even formulate 
correctly in more traditional terms because of the necessity for keeping track 
of boundary conditions, but when recast in the language of $b$-pseudodifferential 
operators, it became much more transparent and amenable to proof, \cite{MP1}. 
There are other, more recent, index-theoretic applications of the
$b$-calculus, such as the higher Atiyah-Patodi-Singer
index theorem 
on Galois coverings by Leichtnam and the second author,
see \cite{LP1}
\cite{LP2}.

Beyond these essentially index-theoretic applications, manifolds
with cylindrical ends, or degenerations to them, have also played
important r\^oles in many other sorts of problems, and by many other
authors; we mention only the important recent work \cite{T}, \cite{MMR},
\cite{KM} \cite{MST} in Donaldson theory and Seiberg-Witten theory.

The two points of view are closely related, and it is by playing them 
off against one another that one may obtain the best insights. To 
illustrate this, we describe in detail a very elementary fact, requisite
to much of what follows, namely, the equivalence between the global 
Atiyah-Patodi-Singer (henceforth APS) boundary conditions for the 
Dirac-type  operator $\eth$ on manifolds with boundary $Z$, assuming 
that all structures are of product-type near the ends, and the natural 
$L^2$ boundary condition on the prolongation $\hZ$ of $Z$ to a manifold 
with infinite cylindrical ends. This equivalence was already noted and used 
in the original paper \cite{APS}.

We begin with a Dirac-type operator $\eth$ over $Z$. To say that $\eth$ is 
of Dirac type means that it acts between sections of two different 
bundles $E$ and $F$, and that there is a parallel bundle map
\[
\gamma: \mbox{\rm Cl\,}(TZ) \longrightarrow \mbox{\rm End\,}(E,F)
\]
which is a fibrewise homomorphism from the Clifford bundle over $Z$ 
to the bundle of endomorphisms from $E$ to $F$, in terms of which, locally, 
\[
\eth = \sum_{j=1}^n \gamma(e_j)\nabla_{e_j} + R,
\]
where $\{e_1,\ldots, e_n\}$ is an orthonormal frame and $R \in
\cin(Z;\mbox{\rm End\,}(E,F))$. Requirements of self-adjointness
impose compatibility conditions on $\gamma$ and $R$. 

Next, suppose that $t$ is a smooth defining function for $Y = \del Z$, 
so that $t$ vanishes simply on $Y$ and is everywhere positive in
the interior of $Z$. A  metric $g$ which is a product near $Y$ 
takes the form 
\[
g = dt^2 + h
\]
where $h$ is a smooth metric on $Y$ which is independent of $t$ for $t 
\leq \e$. The operator $\eth$ is product type near $Y$ if in some collar 
neighbourhood of the boundary, the bundles $E$ and $F$ are lifted
from $Y$, 
\[
E = \pi^* E_Y, \qquad F = \pi^* F_Y, \qquad \mbox{where}\quad
\pi: Y \times [0,t_0) \longrightarrow Y,
\]
and in terms of these partial trivializations, $\eth$ takes the form
\[
\eth = \gamma(\del_t)\left( \del_t + A \right) + R.
\]
Here $A$ is some $t$-independent first order elliptic operator on $Y$,
acting on sections of $E$ and $R$ is is also $t$-independent for
$ t \leq t_0$. In the cases of interest to us, $A$ is self-adjoint,
and we shall also assume that $R \equiv 0$ for simplicity. 
More detailed descriptions of generalized Dirac-type operators are given 
in \cite{tapsit} and \cite{Mue}.

Since $A$ is self-adjoint, elliptic and first order, its spectrum is a 
discrete sequence of real numbers $\{\lambda_j\}$ which is unbounded both 
above and below. The corresponding eigensections of $E$ will be denoted 
$\phi_j$. There is no natural local elliptic boundary condition for this 
operator. The `correct' global boundary condition was one of the very 
important discoveries in \cite{APS}. 
To state it, we first define the orthogonal projection 
\[
\Pi^+_0: L^2(Y;E) \longrightarrow L^2(Y;E)
\]
onto the sum of eigenspaces for $A$ with nonnegative eigenvalues. Thus 
$\Pi^+_0(\phi_j) = 0$ whenever $\lambda_j <  0$ and $\Pi^+_0(\phi_j) = \phi_j$ 
whenever $\lambda_j \geq  0$. The APS boundary conditions involve
letting the operator $\eth$ act on the domain
\begin{equation}
\{u\in\cin(X;E)\; ; \; \Pi^+_0 (u_Y)=0\}.
\label{eq:1.2.1}
\end{equation}
This boundary projection has a classical analogue: 
if $X$ is the disc in ${\CC}$ and $\eth$ is the Cauchy-Riemann operator, then
it is elementary that the restriction of holomorphic functions to
$S^1$ are precisely those with only nonnegative Fourier coefficients.
(Although neither the metric nor the operator are of product type here,
it is not hard to transform them to be of this form.)

To explain this boundary projection better, we consider elements 
$u$ of the nullspace of $\eth$. Solutions of $\eth u = 0$ may be 
analyzed near $Y$ by introducing the eigendecomposition
\[
u(t,y) = \sum_{j} u_j(t) \phi_j(y), \qquad \mbox{ so that } 
\qquad \eth u(t,y) = \gamma(\del_t)\sum_j \left( u_j'(t)
+ \lambda_j u_j(t)\right)\phi_j(y),
\]
valid in the collar neighbourhood ${\cal U} = [0,t_0) \times Y$. Thus if 
$\eth u = 0$, then $u_j(t) = a_j e^{-\lambda_j t}$ for some constants $a_j$
and for all $j$.
Extend the variable $t$, the bundles $E$ and $F$ and the operator $\eth$ to
the manifold
\[
\hZ = Z \cup_Y \left(Y \times (-\infty,0]\right),
\]
obtained by adjoining the half-cylinder ${\RR}^- \times Y$ to $Z$ along the 
common boundary $Y$. Then the solution $u$ automatically extends to
a solution of this equation on $\hZ$. More importantly, $\Pi^+(u(0,\cdot)) = 0$
if and only if this extension decays exponentially. In particular,
elements of the nullspace of the elliptic boundary problem
$(\eth,\Pi^+_0)$ on $Z$ are in one-to-one correspondence with the $L^2$ 
nullspace of $\eth$ on $\hZ$. 

The associated inhomogeneous elliptic boundary problem is 
\[
\eth u = f \quad \mbox{in}\ X, \qquad \Pi^+_0(\left. u \right|_Y) = \phi,
\]
for $f$ and $\phi$ in some appropriate spaces of sections. Assuming that 
$\eth$ is symmetric, the adjoint boundary problem is given by the 
pair $(\eth,\Pi^+)$, where $\Pi^+$ 
is the spectral projection onto the {\it positive}
part of the spectrum of $A$.
Notice that $\ran(\Pi^+_0) \ominus \ran(\Pi^+) =
\ker A$. Elements of the nullspace of $(\eth,\Pi^+)$ 
are in one-to-one correspondence with the {\it extended} $L^2$ nullspace
of $\eth$ on $\hZ$, which contains all temperate solutions
of $\eth u = 0$. To obtain a self-adjoint boundary problem, we
see that we must restrict the domain for this operator to be 
intermediate between the domains for the two boundary problems
$(\eth,\Pi^+_0)$ and $(\eth,\Pi^+)$.
It turns out that the self-adjoint extensions of $(\eth,\Pi^+)$
are in one-to-one correspondence with the Lagrangian subspaces
of $\ker A$. The symplectic structure with respect to which
these subspaces are Lagrangian is the one induced on $\ker A$
from the $L^2$ inner product and the almost complex structure
induced from Clifford multiplication by $\del_t$. (We are using
that $\gamma(\del_t)$ anticommutes with $A$, but in particular
preserves $\ker A$, and that the Clifford relations imply
that $\gamma(\del_t)^2 = -I$.) If $\Lambda \subset \ker A$
is any Lagrangian subspace, then we define the augmented
projection $\Pi^+_\Lambda$ by demanding that $\ran(\Pi^+_\Lambda) \ominus
\ran(\Pi^+) = \Lambda$. The corresponding self-adjoint elliptic 
boundary problem is then given by the pair $(\eth,\Pi^+_\Lambda)$. 

To conclude this discussion, we observe that there is 
a natural Lagrangian subspace $\Lambda_{\mathrm sc} \subset \ker A$
which provides the connection between the `finite'
elliptic boundary problem and the operator on the manifold
$\hZ$. It is the subspace of asymptotic limits of solutions 
of $\eth u = 0$:
\[
\Lambda_{\mathrm sc} \equiv \{ \lim_{t \rightarrow -\infty} u(t,y): 
\mbox{\rm $u$ a bounded solution of $\eth u = 0$ defined
on all of $\hZ$}\},
\]
and is called the {\it scattering Lagrangian} associated
to the operator. 
The fact that the solution $u$ must be globally defined
is very important in this definition. It is not 
immediately clear why this $\Lambda$ should even
have the correct dimension, let alone be Lagrangian, but
these follow from Green's formula and adjointness
considerations, cf. \cite{Mue} \cite{tapsit}. 

\subsection{The calculus of $b$-pseudodifferential operators}

The simple observations of the last subsection indicate that it is 
at least as fruitful to work in the category of manifolds with 
cylindrical ends as of manifolds with boundary, and thus we have reconnected
with our earlier question about the most natural classes of functions, 
operators, etc., on a manifold with boundary. The preceding discussion
points out one advantage to studying geometric
elliptic operators on complete manifolds: there is no need to
choose boundary conditions explicitly, because the $L^2$ requirement 
imposes a natural set of such conditions. 

There are other advantages too. For example, one often underappreciated
fact is that there is a lack of naturality inherent in the usual spaces 
of smooth functions and pseudodifferential operators on a manifold with 
boundary $Z$: general pseudodifferential operator do not 
map the space $\cin(Z)$ to itself.  This prompted Boutet de Monvel's
introduction of the transmission condition for symbols, and his
pseudodifferential calculus adapted to boundary problems in the
mid 1970's \cite{BdM}. Several years later a quite different approach
was initiated by Melrose, resulting in the so-called $b$-calculus, 
or calculus of $b$-pseudodifferential operators on a manifold
with boundary $Z$. In fact, this $b$-calculus is the first
step toward a rather general microlocal approach for studying
a hierarchy of spaces of degenerate differential pseudodifferential
operators. The surgery calculus we discuss later is one
amongst many in this hierarchy, and as can be seen from the geometric 
and analytic problems motivating it, is some sort of extension of 
the $b$-calculus. Because of this relationship, we include a very brief 
introduction to the $b$-calculus here. This will be continued, and the 
surgery calculus itself will be discussed, later in the paper. 

In defining this `calculus' (by which we mean a set of pseudodifferential 
operators which are essentially closed under composition, up to some 
elementary and computable obstructions having to do with integrability of 
certain functions) it is customary to start by introducing the space
of $b$-vector fields on $Z$  (although one might easily regard one of the
other objects we introduce below as the `primary' object of the theory).
This class of vector fields is defined by
\[
{\cal V}_b(Z) = \{V \in {\cal C}^\infty(Z; TZ): V \ \ 
\mbox{\rm tangent to } \ \ \del Z = Y\},
\]
and this condition is obviously closed under Lie bracket, so that
$\Vb(Z)$ is a Lie algebra. In terms of a smooth boundary defining function $x$,
and any choice of local coordinates $y$ on $Y$, $\Vb$ is generated
over $\cin(Z)$ by 
\[
x\del_x, \del_{y_1}, \ldots, \del_{y_k}, \quad k = \mbox{\rm dim\,}Y. 
\]
Alternately, we can also define
\[
\Vb = \{V \in \cin(Z;TZ): Vx = xf\quad \mbox{for some}\ f \in \cin(Z)\}.
\]
Next, a metric $g$ is said to be a $b$-metric if $g(V,W) \in \cin(Z)$,
$g(V,V)\geq 0$ for 
every $V,W \in \Vb$. In terms of the same local coordinates near $Y$,
any such $g$ is a positive smooth symmetric two-tensor in the 1-forms
$dx/x$ and $dy^j$. A slightly more tractable, and just as useful
subclass of these are the {\it exact} $b$-metrics. $g$ is called exact
if there exists
a boundary defining function $x\in\cin(Z)$
and  some smooth (in the ordinary sense) symmetric two-tensor
$h$ such that
\[
g = \frac{dx^2}{x^2} + h.
\]
If we introduce $t = -\log x$, and if $h$ is independent of $x$ in some 
neighbourhood of $Y$, then an exact $b$-metric is nothing more than a 
metric with an infinite product cylindrical end, and a general
exact $b$-metric decreases at an exponential rate to a product
cylindrical metric in these coordinates.

{}From $\Vb(Z)$ we can define the ring of $b$-differential operators 
$\mbox{\rm Diff}_b^{\,*}(Z)$: this contains all operators which may be 
written as locally finite sums of products of elements of $\Vb$. 
Thus, in the same local coordinates, 
\[
\mbox{\rm Diff}_b^{\,m}(Z) \ni L \Longrightarrow L = \sum_{j+|\alpha| \leq m}
a_{j,\alpha}(x,y)(x\del_x)^j\del_y^\alpha.
\]
Examples of such operators  include any geometric operator, such as the 
Laplacian on differential forms or Dirac operator, associated to an 
exact $b$-metric.  

Although any $L \in \mbox{\rm Diff}_b^{\,*}(Z)$ is degenerate in the ordinary 
sense, it is still possible to define a meaningful notion of ellipticity:
a $b$-operator $L$ is said to be elliptic if it may be represented
locally as an elliptic combination of the basic spanning set of
$b$-vector fields listed above. One is then led to ask the 
question as to whether such operators have any `right' to be
called elliptic, i.e. whether they enjoy any of the properties
familiar from elliptic theory on compact manifolds. Specifically,
are these operators Fredholm on any natural function spaces,
and what are the regularity properties of solutions of 
$Lu=0$ or $Lu=f$? Note that because elliptic $b$-operators
are elliptic in the ordinary sense in the interior of $Z$, 
these questions really involve only `local' behaviour at $\del Z$.

To investigate these questions, it is natural to use pseudodifferential
methods, and the heart of this technique is to introduce a class of
pseudodifferential operators $\Psi_b^*(Z)$ which contains 
$\mbox{\rm Diff\,}_b^*(Z)$, and which is also hopefully large
enough to include inverses, or at least good parametrices,
for the elliptic differential $b$-operators. As a clue to
how one might define these pseudodifferential operators, observe
that any element of $\Vb(Z)$, hence any $b$-differential operator,
is approximately invariant under dilations in the variable $x$
(which correspond to translations in the variable $t$). Thus
one might hope to characterize elements of $\Psi_b^*(Z)$ by
this same property, and indeed this is the case. Actually, it
is easiest to characterize these operators by geometric
properties of their Schwartz kernels, which are distributions
on the product $Z \times Z = Z^2$. Because elements of
$\mbox{\rm Diff\,}_b^*(Z)$ are degenerate at $\del Z$, we
expect pseudodifferential operators which represent inverses 
for the elliptic elements to have some sort of singularity at $(\del Z)^2$. 
The main idea is that these singularities can be characterized
geometrically: instead of regarding the Schwartz kernel
of an element $B \in \Psi^*_b(Z)$ as a more
singular distribution on the relatively simple space 
$Z^2$, we instead regard it as a simpler distribution on
a geometrically more complicated space $Z^2_b$, the
$b$-stretched product of $Z$ with itself. This new space
is obtained from $Z^2$ by blowing up the corner $(\del Z)^2$;
said differently, introduce polar coordinates around this
corner and include as part of the blow-up the new
face where the polar distance variable $r = 0$. 
$Z^2_b$ is a manifold with corners. It has three
codimension one boundary faces, $Z \times \del Z$,
$\del Z \times Z$, and the new one created from the
blow-up, and away from this final hypersurface it is
diffeomorphic to $Z^2$. The Schwartz kernel of $B$ is
required to be smooth in the interior of $Z^2_b$ away
from the diagonal, where it is to have an ordinary
pseudodifferential singularity, and at each of the
codimension one boundary hypersurfaces it is required
to have complete polyhomogeneous (i.e. classical) 
expansions. In fact, it is even required to be smooth up to 
the front face. Notice, however, that back on the original
manifold $Z^2$, this means that it is only smooth in
polar coordinates around the corner, away from the diagonal.
When $B \in \mbox{\rm Diff\,}^*_b(Z)$, then its Schwartz
kernel lifts to a $\delta$-section supported along the
lifted diagonal of $Z^2_b$. 

Complete and accessible discussions of this space of operators 
are to be found in \cite{tapsit}, and \cite{Ma}, to which
we refer the interested reader for more details.

To return to the original objective, though, once the space 
 of $b$-pseudodifferential operators on $Z$, $\Psi^*_b(Z)$, 
has been defined, and certain basic facts about it, such
as its closure under composition and a satisfactory symbol 
calculus, have been established, then one may proceed with the investigation 
of the elliptic differential $b$-operators. To phrase the main 
results, one lets these operators act not just 
on the ordinary Sobolev spaces $H^s(Z)$, but rather on 
weighted $b$-Sobolev spaces $x^\delta H^s_b(Z)$. The
subscript $b$ refers to the fact that the differentiations
involved in defining these spaces should be with respect
to the elements of $\Vb$, while the factor $x^\delta$ allows
for changing the rates of growth or decay at $\del Z$.
The basic result is that for all but a discrete 
set of weight parameters $\delta$, an elliptic $b$-operator
$L$ is Fredholm on $x^\delta H^s_b(Z)$. Furthermore, an
arbitrary solution of $Lu = 0$ admits a complete 
polyhomogeneous expansion in powers of $x$ (and possibly
$\log x$) as $x \rightarrow 0$. This polyhomogeneity
is the natural replacement for the special case of
smoothness up to the boundary (which is what occurs
when all powers in the expansion are nonnegative
integers and no logarithmic factors occur). 
These results are all consequences of the fact that
one has a very precise description of a good parametrix
for $L$. In fact, that is really the point of the theory.
After the not completely insignificant effort involved
in defining the calculus and constructing the parametrix,
we then have a more or less complete geometric description of
the Schwartz kernel of the generalized inverse for $L$ on any 
of the admissible weighted Sobolev spaces. From this it is possible
to simply `read off' any more refined mapping or regularity properties
about $L$ one might wish to know. 

This overview of the $b$-calculus is intended to motivate
the similar, but unfortunately more elaborate, description
of the surgery calculus later. 

\subsection{Three different approaches to the problem} 

In this final section of the introduction we briefly introduce 
three different methods which have been developed to study the surgery 
problem. These were all developed roughly simultaneously and independently,
but each was directed toward, and achieved, somewhat different goals. 
What we call the first approach was developed by Bunke \cite{Bu}, the
second by Vishik \cite{V}, and later 
used by Br\"uning and Lesch \cite{BL}, while the final one was contained
in work of the first author with Melrose \cite{asatei}, and then also 
with Hassell \cite{asatet}. To avoid being overly self-referential, 
this last approach will be referred to as that of MM/HMM.  In later 
sections of this paper, we amplify the descriptions of these 
approaches --- rather cursorily for the first 
two and in more detail for the third. For the most part, we shall
only discuss Dirac-type operators because they have provided
the main setting for applications. 

As we have seen, the main issue is to somehow `disconnect' the
operator $\eth_{X_+}$ from $\eth_{X_-}$, and this may be
done either by use of boundary conditions or by literally
disconnecting the two halves geometrically by placing them
at infinite distance from one another.  Bunke's approach
is the least intricate, technically, and essentially uses 
both of these types of considerations. Vishik's ideas involve
a variation of boundary conditions along $H$, while
those of MM/HMM rely on the idea of geometric separation.

The goals of these papers are also quite different. Bunke's
intent is to find a gluing formula for the eta invariant. Vishik was 
concerned with gluing formul\ae\ for determinants of elliptic operators, 
and particularly for the analytic torsion. In the somewhat more tractable 
version of his ideas developed by Br\"uning
and Lesch, the goal is to find another proof of the gluing formula for 
the eta invariant.  The techniques of MM/HMM are directed toward 
proving uniformity of the resolvent and heat kernel associated to $\eth^2$ 
in the `analytic surgery limit' as the manifold $X$ stretches to infinite 
length along the hypersurface $H$. Gluing formul\ae\ for the eta invariant, 
and analytic torsion \cite{H}, are then consequences of this uniformity, but 
far more detailed information is obtained along the way. The expense, of 
course, is that the development of this approach is by its nature the most 
technically intricate of the three.

We now go into only slightly more detail. Bunke's setup involves
considering two different manifolds. The first is the disjoint
union $X_+ \sqcup X_-$, each endowed with long (but finite)
cylindrical ends, while the second is the disjoint union
of $X$, endowed with a long cylindrical section around $H$,
and a long cylindrical piece $[-T,T] \times H$. The goal is
to show that there is some abstract unitary equivalence between the 
Dirac operators on these two (sets of) manifolds; since they
are unitary invariants, the eta invariants for these manifolds
must also coincide. Three of the four components are the various 
terms one expects in the gluing formula for the eta invariant, and the
fourth represents the defect term, and it may be computed
`explicitly'. One subtlety here is in determining how the 
different boundary conditions at the various ends arise in 
order that the unitary isomorphism be valid. 

Vishik's setup, on the other hand, involves the consideration
of a family of boundary conditions along the hypersurface $H$.
Each corresponds to some self-adjoint elliptic boundary
problem. At one extreme, this problem corresponds to the
operator $\eth$ on the closed manifold $X$, where the 
hypersurface $H$ becomes `invisible'; these are the transmission
boundary conditions. At the other extreme, the boundary conditions
are the natural APS ones on each half. The analytic torsion,
or eta invariant, may be computed for each operator in this
family, and the problem then consists of computing the
variation of these invariants with respect to the parameter.
The total variation with respect to this parameter represents
the defect. 

Finally, as indicated earlier, in MM/HMM the goal is to develop
a `surgery calculus' $\Psi_s^*(X)$, that is, a calculus of 
pseudodifferential operators on $X$, depending on a parameter $\e$, 
and which incorporates the sorts of degeneracies seen in the family
of Dirac operators or Laplacians with respect to metrics
undergoing degeneration to infinite cylindrical ends. Thus,
for $\e > 0$, the surgery calculus restricts to the ordinary 
pseudodifferential calculus, while at $\e = 0$ it somehow
induces the $b$-calculus. The point is to show how the transition 
between these quite different calculi takes place. Once 
this calculus is defined, and its basic analytic properties established, 
such as a symbol calculus, closure under composition, etc., then one may 
use it to construct parametrices for $(\eth_{X,g_\e}^2-\lambda)$, for
example.  As with the $b$-calculus, if one is able to describe the 
behaviour of this resolvent, or of the heat kernel, uniformly 
with respect to $\e$, and explicitly, it is then straightforward to 
examine the behaviour of these auxiliary numerical invariants. 
One also obtains more detailed information, such as the
way in which the discrete spectrum accumulates into
continuous spectrum. 
 
\section{Some applications of the surgery formula -- an overview}  

We now describe in somewhat greater detail four different types
of objects, for the study of which some form of the analytic
surgery technique has proved useful. These are index bundles,
the eta invariant, analytic torsion and determinant bundles. 

In each of the following settings we shall, again for simplicity,
consider only the case of a Dirac-type operator $\eth$, acting
between sections of the bundles $E$ and $F$ over the manifold
$X$. We shall describe at least the general form of the gluing
theorems in each context, leaving the more precise statements
until later. 
 
\subsection{Index bundles}  
First we consider the numerical index. We assume that $X$ is 
{\it even}-dimensional and, for simplicity, spin. We denote by 
$\Spi=\Spi^+\oplus \Spi^-$ the spin bundle and its splitting
into the plus and minus spin bundles. The Dirac operator $\eth$ is 
odd with respect to the natural ${\ZZ}_2$-grading, and so takes the
form 
\[ 
\left( 
\begin{array}{ccc} 
 0 & \eth^- \\ 
        \eth^+ & 0 \end{array} \right) 
\quad\quad \eth^-=(\eth^+)^* 
\] 
with $\eth^\pm: \cin(X,{\Spi}^{\,\pm})\rightarrow \cin(X,{\Spi}^{\,\mp})$.  
Since the manifold $X$ is closed and compact the Dirac operator 
is Fredholm on any Sobolev space $\eth^+: H^m(X,\Spi^{\,+})\rightarrow 
H^{m-1}(X,\Spi^{\,-})$ with index $\ind(\eth)=\dim(\ker(\eth^+)) 
-\dim(\ker(\eth^-))$ which is independent of $m$ by elliptic regularity. 
In this even-dimensional setting, our primary interest 
is in $\eth^+$ and not in the full self-adjoint operator $\eth$.  
 
Suppose now that $X=X_+\cup_H X_-$; up to a perturbation not affecting 
the index we can assume that the metric near the disconnecting 
hypersurface $H$ is of product type. This means that the Dirac 
operator on each piece $X_{\pm}$ takes the product form 
introduced in the previous section, with $t$ now denoting a 
defining function for $H$. Notice that the vector field $\del_t$ will 
be normal to $H$, the common  boundary of $X_\pm$, but inward pointing 
for one manifold, say $X_+$, and outward pointing for the other. 
As in \S 1.2, we denote by $\Pi^+_0$ the augmented APS spectral 
projection for the boundary operator $\eth_H$. Because of the discrepancy 
in the orientation of the normals it is easy to check directly
that the two APS boundary value problems can be written as 
$(\eth_{X_+}^+,\Pi^+_0)$ and $(\eth_{X_-}^+,\mbox{\rm Id}-\Pi^+)$. 
Applying the Atiyah-Singer index theorem
\cite{AS3} to $\eth_X$ and the Atiyah-Patodi-Singer 
index theorem \cite{APS} to the two boundary value problems and observing 
moreover that the two eta invariants cancel because of the opposite 
orientation of the normals, we obtain the following surgery formula for 
the index 
\begin{equation} 
\ind(\eth_X)=\ind(\eth_{X_+},\Pi^+_0)+\ind(\eth_{X_-}, 
\mbox{\rm Id}-\Pi^+)+\dim(\mbox{\rm null\,}\eth_H). 
\label{eq:2.1.1} 
\end{equation} 
Notice that $\dim(\mbox{\rm null\,}\eth_H)=\dim(\Pi^+_0-\Pi^+)$. 
Suppose now that $X$ is the typical fibre of a fibration 
of compact manifolds: $\phi: M\rightarrow B$. Following what
is now standard notation, we denote the family of metrics on the 
fibres by $g_{M/B}$; we also assume that the fibres carry smoothly 
varying spin structures.  We denote by $M^z$ the fiber over 
$z$; thus  $M^z\equiv \phi^{-1}(z)\cong X$. For each $z\in B$ we can   
consider the Dirac operator $(\eth_{M})_z$ naturally defined by the spin 
structure of $M^z$. We obtain  a family $\eth_M=((\eth_M)_z)_{z\in B}$  
of Dirac operators. 
 
Let  $H$ be a disconnecting hypersurface of the fibration $M\rightarrow 
B$ and assume that $H$ also fibres over $B$: $\phi|_{H}: H\rightarrow B$. 
Thus $M=M_+\cup_{H} M_-$ and each fiber $M^z$ is  the union along $H^z$ 
of two manifolds with boundary: $M^z_+\cup_{H^z} M_-^z$. 
 
We obtain in this way four families of Dirac operators: 
$\eth_M$, $\eth_{M_+}$, $\eth_{M_-}$ and $\eth_H$. 
Since the family $\eth_M$ is defined on a closed manifold,
the familiar construction of an {\it index bundle},
as in \cite{AS4}, provides us with  an index class $\Ind(\eth_M)\in K^0(B)$.
 The problem is to formulate and 
prove the analogue of (\ref{eq:2.1.1}). If the vector spaces 
$\ker(\eth_{H})_z$ are of {\it constant rank} as $z$ varies in $B$, 
then the APS boundary value problems for $M_+^z$ and $M_-^z$ define, 
as $z$ varies, two {\it continuous} families of boundary value problem 
\[ 
(\eth_{M_+},\Pi^+_0)=(\eth_{M^z_+},(\Pi^+_{0})_z)_{z\in B} 
\quad\quad 
(\eth_{M_-},\Id-\Pi^+)=(\eth_{M^z_-},\Id-(\Pi^+)_{z})_{z\in B}. 
\] 
Notice that because of the assumption of constant rank, $\Ker(\eth_H)= 
\cup_{z\in B} \ker(\eth_{H})_z$ is a {\it continuous} (in fact smooth) 
vector bundle. 
 
However if the constant rank assumption is not satisfied we know  
that the two APS families will not be continuous. In this case, 
as explained at length in Part I of this survey \cite{P1}
we need to fix a {\it spectral section} $P$ for the family $\eth_H$ 
(see \cite{MP1} \cite{MP2}). Thus $P=(P_z)_{z\in B}$ is a smooth family of   
pseudodifferential operators of order zero which are self-adjoint projections 
and with the additional property that there exists a positive constant 
$R\in \RR$ such that 
\begin{equation} 
(\eth_{H^z}) u=\lambda u\Rightarrow  
\left\{\begin{array}{ll} 
P_z u=u & \mbox{if $\lambda>R$}\\ 
P_z u=0 & \mbox{if $\lambda<-R$} 
\end{array} \right. 
\label{eq:2.1.3} 
\end{equation} 
Since $\eth_H$ is by construction a boundary family we know  
from \cite{MP1} that there exist an infinite number of spectral sections;  
moreover two  spectral sections, $P_1$, $P_2$, give rise to a difference 
element $[P_1-P_2]\in K^0(B)$ (and in fact it can be proved that 
these differences exhaust all of $K^0(B)$). A spectral section $P$ for 
$\eth_H$ fixes a {\it smooth} family of generalized APS boundary value problem 
$(\eth_{M_+},P)$ and thus an index class $\Ind(\eth_{M_+},P)$. 
Simply define $(\eth_{M_+},P)_z$ as the operator $\eth_{M_+^z}$ with domain 
\[ 
\{u\in L^2(M^z_+,E_z) ; (\eth_{M^z_+}) u\in L^2(M^z_+,E_z),  
P(u|_{\del M^z_+})=0\}. 
\] 
We also obtain an index class for $\eth_{M_-}$ by considering the family 
of boundary value problems $(\eth_{M_-},\Id-P)$ (recall that the normal to 
$M_-$ is oriented in the outward direction). 
 
We can now state the decomposition formula for index bundles (\cite{DZ}).  
 
\begin{theorem} (Dai-Zhang) 
Let $P_1$ and $P_2$ be spectral sections for $\eth_H$. 
Then the following formula holds 
\[ 
\Ind(\eth_M)=\Ind(\eth_{M_+},P_1)+\Ind(\eth_{M_-},\Id-P_2) 
+[P_1-P_2]\quad\mbox{\rm in}\quad K^0(B) 
\] 
\end{theorem} 
 
The corresponding formula for the Chern characters follows directly from the 
family index formula of Melrose-Piazza. Notice that if in particular 
$\Ker(\eth_{H})$ is a smooth vector bundle then we can choose $P_1=\Pi^+_0$, 
$P_2=\Pi^+$; then $[\nul(\eth_H)]=[P_1 - P_2]$ and we  obtain the precise 
analogue of (\ref{eq:2.1.1}). 
 
The proof of the surgery formula for the index bundle, as given by 
Dai-Zhang, follows  Bunke's method. One can also prove it with
the surgery calculus of MM, as in \cite{P3}.
 
\subsection{The eta invariant}  
The eta invariant $\eta(\eth)$ is a spectral invariant which was discovered 
originally in the context of the APS index theorem as the boundary correction 
term. Because indices are really even-dimensional phenomena, 
eta invariants are therefore of most interest when the dimension  
of $X$ is odd -- in fact, there are a number of difficult 
analytic subtleties in even defining the eta invariant 
in even dimensions. It was pointed out by Singer \cite{S1} 
that the eta invariant of Dirac operators in odd dimensions 
actually shares many properties with the index of Dirac 
operators in even dimensions: it could even be regarded 
as the odd-dimensional analogue of the index. This has  
been made more precise and generalized considerably by  
Melrose \cite{Me3}. Other aspects of Singer's assertions 
have been addressed and proved by Wojciechowski \cite{Woj1} 
\cite{Woj2}. 
 
In any case, the basic definition is given as follows. Since  
$\eth$ is elliptic, and is self-adjoint when $n \equiv \dim X$ is odd,  
its spectrum is discrete, and we denote it by $\{\lambda_j\}$. 
This sequence is unbounded in both directions. The Weyl asymptotics 
show that the eta function 
\[ 
\eta(s) \equiv \sum \frac{\mbox{\rm sgn\,} \lambda_j}{|\lambda_j|^s} 
\] 
is defined and holomorphic in the half-plane $\mbox{\rm Re\,}s > n$.  
From the usual analysis of the short-time asymptotics of the heat
kernel associated to $\eth^2$, this function extends meromorphically 
over the complex plane. This extension is actually regular at $s=0$ 
(and this is precisely the point that becomes much more subtle for 
non-Dirac operators and in even dimensions). This value is defined to 
be the eta invariant $\eta(\eth)$.  
 
It is usually easier to work with a different expression for this invariant. 
First note that for $\mbox{\rm Re\,}s$ sufficiently large,  
\[ 
\eta(s) \equiv \frac{1}{\Gamma(\frac{s+1}{2})} 
\int_0^\infty 
t^{\frac{s-1}{2}}\mbox{\rm Tr\,}\left(\eth e^{-t\eth^2}\right)\,dt. 
\] 
The integral continues to converge up to $s=0$. This is 
not obvious, but follows from Getzler's rescaling technique, 
which is explained in \cite{BGV} and \cite{tapsit}. The  
factor in front is regular at $s=0$ too, and so  
\[ 
\eta(\eth) = \frac{1}{\sqrt{\pi}}\int_0^\infty t^{-1/2}\mbox{\rm Tr\,} 
\left(\eth e^{-t\eth^2}\right)\,dt. 
\] 
 
Since we are also interested in the eta invariant on the 
manifolds with boundary $X_\pm$, we must also discuss how to 
define these. Recall the two different sorts of metrics 
we have been considering, those which make these manifolds 
compact with boundary, and are of product type near 
the boundaries, and those which are complete with infinite 
cylindrical ends. In the former case, following our earlier 
discussion, once we have introduced Lagrangian  
subspaces $\Lambda_\pm \subset \ker \eth_{H}$, we obtain 
self-adjoint elliptic boundary problems $(\eth_{X_+},\Pi^+_{\Lambda_+})$ 
and $(\eth_{X_-},\Pi^+_{\Lambda_-})$. These operators 
have discrete spectrum, and at least formally the preceding 
definitions make sense. The details of making these plausible 
definitions rigorous has been carried  
by M\"uller \cite{Mue}. The invariants we obtain will 
be denoted $\eta(\eth_{X_{\pm}},\Lambda_\pm)$ for simplicity.  
 
We may now state one form of the surgery problem for eta invariants  
explicitly: find a tractable expression for the defect  
\[ 
\delta(\Lambda_+,\Lambda_-) \equiv \eta(\eth_X) -  
\eta(\eth_{X_+},\Lambda_+) - \eta(\eth_{X_-},\Lambda_-). 
\] 
Notice that we have written this defect as a function of 
the two Lagrangians $\Lambda_\pm$. That these should be 
the essential variables on which it depends requires 
some work. Furthermore, it is also of substantial interest 
to see whether there is some choice of Lagrangians 
for which this formula becomes particularly simple or 
natural.   
 
The other geometric setup is when $X_\pm$ are endowed 
with exact $b$-metrics, and the degeneration from $X$ to  
these complete metrics is a surgery degeneration as 
defined more precisely later in this paper. The biggest 
obstacle now is that the eta invariant $\eta(\eth)$  
ostensibly requires the operator $\eth$ to have  
discrete spectrum, while in this case we know that the 
Dirac operator has continuous spectrum. Thus it is not 
even clear how to define the eta function $\eta(s)$. 
Starting from the definition in terms of the heat 
kernel we see the obstacle in a different way.  
The heat kernel for $\eth^2$ on either $X_+$ or $X_-$ is a  
smooth function, but it does not decay along the 
cylindrical ends, and hence is not integrable.  
The way out of this impasse is to use the regularized 
$b$-trace defined by Melrose \cite{tapsit}. This $b$-trace 
is an extension of the ordinary trace to the ring 
of smoothing $b$-pseudodifferential operators. For 
any such operator $R$, the pointwise trace along 
the diagonal of the Schwartz kernel of $R$, (i.e. just 
its restriction to the diagonal, if $R$ acts 
on functions) has an asymptotic expansion in terms 
of nonnegative powers of the boundary defining function $x$. 
It therefore makes sense to define 
\begin{equation}
{}^b\mbox{Tr\,}(R) = \lim_{\e \rightarrow 0} 
\left(\int_{x \geq \e} K_R(x,y)\,\frac{dxdy}{x} + 
\log\e\cdot\int_{x=0} K_R(0,y)\,dy\right) 
\label{eq:2.2.1}
\end{equation}
Here $K_R(x,y)$ is the (pointwise trace of the) 
Schwartz kernel of $R$ on the diagonal, using some local 
coordinate system $y$ on the boundary. Granting the 
naturality of this definition, it is then reasonable 
to define the regularized $b$-eta invariant, ${}^b\eta(\eth)$, 
via the same heat-kernel formula, but substituting 
the $b$-trace for the ordinary trace.  
 
Although this definition may seem ad hoc, it follows from work of 
M\"uller \cite{Mue} that on a manifold $\hZ$ with infinite (product) 
cylindrical ends, if $\Lambda_{\mathrm sc}$ is the scattering 
Lagrangian, i.e. the subspace of $\ker \eth_Y$ obtained as the set of 
asymptotic limits of $\eth u = 0$ on $\hZ$, and if $Z_T$ represents 
the truncation of $\hZ$ to any compact piece, where $t \geq T$ for 
some sufficiently negative $T$,
 then 
\[ 
\eta(\eth_{Z_T},\Pi^+_{\Lambda_{\mathrm sc}}) = {}^b\eta(\eth_{\hZ}). 
\] 
In particular, the term on the left does not depend on the length of the 
cylindrical ends. This serves as ample evidence that the $b$-eta
invariant is a natural object.
 
The easiest case to calculate this defect is when $\ker \eth_H  
= \{0\}$. Then necessarily both $\Lambda_+$ and $\Lambda_-$ 
are also trivial, and so it is hardly surprising (and is 
consonant with our notation) that the defect vanishes. 
In this setting of boundary problems on finite manifolds this  
was proved by Bunke \cite{Bu}, while in the setting of  
degeneration to exact $b$-metrics it was proved in \cite{asatei}.  
 
It is much more interesting, of course, to see what happens when 
$\ker \eth_H$ is nontrivial. Then $\delta(\Lambda_+,\Lambda_-)$ 
does not necessarily vanish. Rather nicely, it turns out 
that the `best' Lagrangians with respect to which to compute this 
defect are the ones induced from the asymptotic limits 
of bounded solutions on the cylindrical extensions of  
$X_\pm$. The importance of these Lagrangians even for 
the boundary problem on the compact manifold was first 
proved by M\"uller \cite{Mue}. There are several nice 
formul\ae\ for the defect in this case. The one obtained 
by Bunke is in the form of an averaged Maslov invariant, 
\[ 
\delta(\Lambda_+,\Lambda_-) = \int_G \mu(g\Lambda, \Lambda_+, 
\Lambda_-)\,dh. 
\] 
Here $G$ is the Lagrangian Grassmanian of $\ker \eth_H$, 
$dh$ is normalized Haar measure on it, $\mu$ is the  
Maslov invariant, which is a function of three separate 
Lagrangians, and $\Lambda$ is an arbitrary third Lagrangian. 
The formula found for the defect in this case in \cite{asatet} is  
more elementary and purely linear algebraic. Again 
when $\Lambda_\pm$ are the scattering Lagrangians, 
\begin{equation}
\delta(\Lambda_+,\Lambda_-) = \frac{i}{2\pi} \log \mbox{\rm Sdet\,}
(I - S_+^r S_-^r).
\label{eq:elexp}
\end{equation}
Here 
$\mbox{\rm Sdet\,}$ is the superdeterminant
of a ${\ZZ}_2$-graded diagonal operator,
\[
\mbox{\rm Sdet\,} 
\left( 
\begin{array}{ccc} 
 A^0 & 0\\ 
        0 & A^1 \end{array} \right) 
\, \equiv \det A^0\cdot (\det A^1)^{-1} \,\,
\] 
and 
$S_{\pm}^r$ are differences of certain projection 
operators associated to the scattering Lagrangian
subspaces $\Lambda_{\mathrm sc}$ in $\mbox{ker\,}(\eth_H)$.
We refer to \cite{asatet} for further explication. 

A similar, but less tidy, finite-dimensional linear
algebraic formula for this term was discovered (earlier) 
by
Lesch and Wojciechowski \cite{LW} in their study of the
closely related problem of the dependence of the eta invariant 
on cylinders with boundary conditions given by arbitrary
Lagrangian subspaces. It arises, as does the expression
above, from one other useful way to think of this defect, 
namely as the eta invariant of an associated one-dimensional 
problem. The Dirac-type operator $\gamma \del_s$ acts on the 
space of sections of the trivial bundle $I \times 
\mbox{\rm ker\,}(\eth_H)$ over $I = [-1,1]$, where $\gamma$ 
denotes Clifford multiplication by the unit normal to $H$. 
Imposing the boundary conditions associated with the
Lagrangian subspaces $\Lambda_{\pm}$ at $s = \pm 1$, 
we obtain a self-adjoint operator with discrete spectrum.
Then it is also true that the defect is simply the
eta invariant of this operator, 
\begin{equation}
\delta(\Lambda_+,\Lambda_-) = \eta(\gamma \del_s, \Lambda_+,
\Lambda_-).
\label{eq:onedeta}
\end{equation}
The seemingly more explicit expression (\ref{eq:elexp})
is deduced from this one.
 
\subsection{Analytic torsion}  
 
The analytic torsion was first defined in the seminal 
paper of Ray and Singer \cite{RS}. In order
to introduce it we first need
to recall the definition of zeta function 
$\zeta_P(s)$ associated to a (second order) self-adjoint elliptic 
differential operator $P$ on a closed manifold $Z$ \cite{S0}.
  This zeta function is defined in manner 
similar to, but simpler than, the eta function. Thus, 
if $\mbox{\rm spec\,}(P) = \{\lambda_j\}$, now a sequence 
of real numbers tending only to infinity, then 
\[ 
\zeta_P(s) \equiv \sum_{j=0}^\infty \lambda_j^{-s}. 
\] 
We are assuming here that all eigenvalues are positive. 
If there are finitely many nonpositive ones we simply 
omit them from this sum. The definition that Ray and 
Singer gave to the determinant is 
\[ 
{\det}' P = e^{-\zeta_P'(0)}. 
\] 
The notation ${\det}'$ is meant to indicate that the 
nonpositive eigenvalues (in particular, the zero 
eigenvalue) have been omitted. 
 
This determinant has emerged as one of the central 
objects of study in spectral geometry (and also 
plays a prominent r\^ole in string theory). In particular, 
it was a crucial ingredient in the proof of compactness 
of isospectral sets by Osgood, Phillips and Sarnak,
and many new results about it have been obtained and 
applied by many authors over the past decade. We mention 
only the work of S.-Y. Chang, Yang, Gursky, Okikiolu, in addition 
to that of Vishik. 
 
Just like the eta invariant, ${\det}' P$ depends on the 
Riemannian metric $g$: 
its variation with respect to a family of metrics is actually  
computable as an integral of local quantities. Ray and  
Singer observed, though, that if one computes the determinants 
$\det'(\Delta_p)$ of the Hodge-Laplacian on $p$-forms, 
and takes a certain weighted sum of these expressions, 
then the resulting object has much more invariance. This 
weighted sum is the analytic torsion $T(Z,g)$ 
defined by 
\[ 
\log T(Z,g) 
= \sum_{p=0}^n (-1)^p p \, {\det}'(\Delta_p). 
\] 
This expression may also be defined when the Hodge-Laplacian 
is twisted by some flat bundle $E$, and this yields 
a number $T(Z,E,g)$. Actually, this expression is independent  
of the metric only when the twisted de Rham complex is 
acyclic, i.e. has all cohomology groups vanishing. Notice 
that this never occurs in the untwisted case. In the 
general case there is an explicit factor which contains 
all metric dependence.   

To be somewhat more precise fix a set of bases $\{\mu\} = 
\{\mu^{(i)}_j\}$ for the cohomology spaces $H^i(Z)$. Next, 
using the Hodge theorem, let $\{\omega^{(i)}_j\}$ be a basis of 
harmonic forms for each of 
these spaces which is orthonormal with respect to the $L^2$-inner 
product induced by $g$. Now let $\Lambda(g,\{\mu\})$ be the 
determinant of the change of basis matrix. Following \cite{RS} 
we define $T(M,\{\mu\})$ as the product $T(Z,\{\mu\})\equiv 
T(Z,g)\cdot \Lambda(g,\{\mu\})$. It is this quantity which is {\it independent}
of the metric $g$. Because of this somewhat surprising invariance, 
and for other (more compelling) reasons, Ray and Singer  
made their famous conjecture that the analytic torsion
$T(Z,\{\mu\})$ agrees with the Reidemeister torsion $\tau(Z,\{\mu\})$, 
a PL invariant of the manifold $Z$ defined by Reidemeister and Franz in the 1930's.  
This conjecture was proved in the late 1970's independently  
by Cheeger \cite{Ch0} and M\"uller \cite{Mue0}.  

By now there are numerous proofs of this Cheeger-M\"uller 
theorem. Many of these rely on ideas 
related to surgery degeneration, including Cheeger's original proof
\cite{Ch0}, and also the fairly recent proof by Burghelea, Friedlander
and Kappeler, which is based on Witten's deformation method
\cite{BFK}. M\"uller's proof \cite{Mue0} was of a somewhat different
nature.  The proofs we highlight here are those by Hassell \cite{H}, 
who applied the surgery calculus of MM/HMM to obtain a gluing formula for 
the analytic torsion, and hence could then follow Cheeger's strategy
to prove its equivalence with Reidemeister torsion, and Vishik's 
gluing formul\ae\  for determinants \cite{V}, the proofs of which employ his 
method of changing boundary conditions, as described below in the slightly
different context of gluing formul\ae\ for eta invariants (developed
by Br\"uning and Lesch \cite{BL})). 

It is obviously not surprising that the results and methods for obtaining 
gluing formul\ae\ for determinants are much the same as for eta invariants.  
Both are obtained from integrating the appropriate heat kernel, or combinations
thereof.  At any rate, the studies of these two quantities are intimately 
interrelated. For reasons of space and time, we shall concentrate almost
exclusively on results obtained for the eta invariant in the remainder of 
the paper. 
 
\subsection{Determinant bundles}    
 
Let $E$ be a finite dimensional vector space and suppose that $T:E\rightarrow E$ 
is linear. Then $T$ induces in a natural way a map $\det(T):\Lam(E)
\rightarrow \Lam (E)$ and hence an element $\det(T)\in (\Lam(E))^*\otimes 
\Lam(E)$; the numerical determinant of $T$ is simply obtained by fixing a 
basis of $E$, or at least a nonzero element of $\Lam(E)$. If $T\in \Hom(E,F)$, 
with $\dim E= \dim F$,  then $\det(T)$ is again well defined as an element of 
$(\Lam(E))^*\otimes \Lam(F)$. 
 
The natural exact sequence of vector spaces 
\[ 
0\rightarrow\ker T\rightarrow E\rightarrow F\rightarrow\coker T\rightarrow 0 
\] 
induces a natural isomorphism 
\[ 
(\Lam(E))^*\otimes \Lam(F)\cong (\Lam \ker T)^*\otimes (\Lam \coker T). 
\] 
Suppose now that the vector spaces $E, F$ and the linear map 
$T$ depend smoothly on a parameter $z\in B$. In other words suppose 
that $T\in\cin(B,\Hom(E,F))$. Applying the preceeding 
remarks we obtain a smooth section $\det T\in\cin(B,{\cal L})$  
of the {\it determinant line bundle} 
${\cal L}$ with fiber at $z$ equal to  
${\cal L}_z=(\Lam(E_z))^*\otimes \Lam(F_z)$. 
 
Notice again that for each fixed $z\in B$ 
\begin{equation} 
{\cal L}_z\cong (\Lam \ker T_z)^*\otimes (\Lam \coker T_z). 
\label{eq:2.5.1} 
\end{equation} 
These elementary remarks show that the determinant of a family 
of linear maps is only defined as a {\it section} of a complex line bundle. 
Of course if would be desirable to have a determinant 
{\it function} $\mbox{\rm DET}: B\rightarrow \CC$ 
assigning a number to each linear map $T_z$. 
If ${\cal L}$ is trivial this can certainly be done 
by fixing a trivializing section $\tau\in\cin(B,{\cal L})$ 
and comparing $\det T$ and $\tau$, viz. 
\[ 
(\det T)(z)=\mbox{\rm DET}_\tau (T_z) \tau(z). 
\] 
The determinant function $\mbox{\rm DET}_\tau$ so obtained depends of 
course on the trivializing section $\tau$ and it is natural to ask 
whether it possible to agree on a {\it canonical} choice. This is 
particularly important in physics. One way to proceed would be to assume 
that the determinant line bundle is equipped with a metric and 
compatible connection. Using a nontrivial covariant constant section 
$\overline{\tau}$ to trivialize this bundle would fix the determinant 
function up to a global phase $C\in U(1)$. The local and global
obstructions to the existence of $\overline{\tau}$ are given by the 
curvature and holonomy of the given connection, respectively. In the
usual physics parlance these are called the local and global anomalies.
 
Again partly motivated by physics, the problem arises as to whether
these ideas can be extended to the infinite dimensional context where
the operator $T$ is replaced by a family of Dirac operators.
Thus let $\eth=(\eth_z)_{z\in B}$ be a smooth family of Dirac operators, 
associated as in \S 2.1 to a smooth fibration of closed compact 
manifolds $\phi:M\rightarrow B$ with even dimensional fibres. 
Since each $\eth_z$ is Fredholm, it makes sense to define the complex line  
\[ 
\Lam\ker(\eth_z)^* \otimes \Lam\coker(\eth_z). 
\] 
However, since the kernel and cokernel of $\eth_z$ are not constant in 
$z$, these complex lines do {\it not} vary smoothly with the parameter $z$. 
The first step then is to show that there exists a smooth complex line bundle 
${\cal L}(\eth)$ over $B$ with the property that its fiber over $z\in B$ is 
naturally identified with $\Lam\ker(\eth_z)^* \otimes \Lam\coker(\eth_z)$. 
The second step is to introduce a metric and compatible connection on 
${\cal L}(\eth)$ in some natural way, and then to compute in geometric terms 
the curvature and the holonomy, i.e. the local and global anomaly. 

This program was accomplished by Quillen in his seminal paper \cite{Q} 
in the special case of $\overline{\del}$-operators on Riemann surfaces 
acting on a vector bundle $E$ and with parameter space ${\cal A}$ equal to  
the moduli space of holomorphic structures on $E$. Since ${\cal A}$ is 
simply connected, only information about the curvature of the so-called 
Quillen metric is required to determine whether the determinant line bundle 
may be trivialized by parallel transport. These results of Quillen were 
extended to the general case by Bismut and Freed in two papers \cite{BF}
(see also \cite{Ch2}, \cite{Wi}, \cite{Fr}, \cite{DaiF}). 
We now illustrate a few of the main ideas behind these works. 

First we define the determinant line bundle ${\cal L}(\eth)$ associated 
to the family $\eth=(\eth_z)_{z\in B}$. We only treat the Dirac case here,
but this construction can be applied to any family of Fredholm operators.
 
Since the fibers  of $\phi:M\rightarrow B$ are even dimensional, each 
Dirac operator may be written as 
\[ 
\left( 
\begin{array}{ccc} 
 0 & \eth_z^- \\ 
        \eth_z^+ & 0 \end{array} \right) 
\quad\quad \eth_z^-=(\eth_z^+)^* 
\] 
with $\eth_z^\pm: \cin(M^z,\Spi^{\,\pm}_z)\rightarrow \cin(M^z,
\Spi^{\,\mp}_z)$. If $E=E^+ \oplus E^-$ is a ${\ZZ}_{2}$ graded 
vector space we use the notation $\det(E)$ for the complex line 
$\Lam(E^+)^*\otimes \Lam(E^-)$. 
 
Clearly if  $(\ker(\eth^\pm_z))_{z\in B}$ form two smooth vector 
bundles, $\Ker(\eth^+), \Ker(\eth^-)$, then the line bundle 
\[ 
{\cal L}(\eth)=\Lam \Ker(\eth^+)\otimes \Lam \Ker(\eth^-)=\det(\Ker \eth) 
\] 
is globally well defined. Notice that if $\Lap^\pm_z=\eth^\mp_z \eth^\pm_z$ 
then it is also true that ${\cal L}({\eth})=\det(\Ker \Lap)$. 
In general consider the  set $U_\lambda=\{z\in B ; \lambda\notin 
\spec(\Lap_z)\}$. Since the spectrum of each Laplacian is discrete, this is 
either a non-empty open set, or else the empty set. Since the latter
may happen for at most a countable set of values of $\lambda$, we may
cover $B$ by a finite collection of such sets $U_{\lambda_k}$.  
Let $\Pi^\pm_{[0,\lambda)}(z)$ be the spectral projection associated to 
the interval $[0,\lambda)$ for the Laplacian $\Lap^\pm_z$. Consider 
\[ 
H^\pm_{[0,\lambda)}(z)=\Ima P^\pm_{[0,\lambda)}(z). 
\] 
This is simply the direct sum of the eigenspaces of $\Lap^\pm_z$
associated to the eigenvalues in $[0,\lambda)$.
As $z$ varies in $U_\lambda$ these vector spaces define a ${\ZZ}_2$-graded
smooth vector bundle $H_{[0,\lambda)}= H^+_{[0,\lambda)}\oplus H^-_{[0,
\lambda)}$. We {\it define} ${\cal L}(\eth)$ restricted to $U_\lambda$ as 
$\det(H_{[0,\lambda)})$. This is a smooth complex line bundle over $U_\lambda$ 
and moreover for each fixed $z\in U_\lambda$ there is a natural isomorphism 
\[ 
({\cal L}(\eth))_z \equiv \det(H_{[0,\infty)}(z))\cong (\Lam (\ker \eth_z^+))^* 
\otimes \Lam (\ker \eth^-_z). 
\] 
coming from the exact sequence 
\[ 
0\rightarrow\ker\Lap^+_z\rightarrow H^+_{[0,\lambda)}(z) 
\rightarrow H^-_{[0,\lambda)}(z)\rightarrow \ker\Lap^-_z\rightarrow 0. 
\] 
Now if $\mu>\lambda$, then on $U_\lambda \cap U_\mu$ we have $H_{[0,\mu)}=
H_{[0,\lambda)}\oplus H_{[\lambda,\mu)}$ and thus $\det(H_{[0,\mu)})\cong
\det(H_{[0,\lambda)}) \otimes\det(H_{[\lambda,\mu)})$. Moreover the 
restriction  of $\eth^+_z$ to $H^+_{[\lambda,\mu)}(z)$ is an isomorphism 
for each $z\in U_\lambda \cap U_\mu$; this means that we can identify 
$\det(H_{[0,\mu)})$ and $\det(H_{[0,\lambda)})$ over $ U_\lambda \cap U_\mu$ 
using the non-vanishing section $\det((\eth^+)_{[\lambda,\mu)})$. 
The resulting  line bundle, which is now defined over all of $B$ is, by 
definition, the determinant line bundle ${\cal L}(\eth)$ defined by the 
family $\eth$. By construction there is a natural isomorphism 
\[ 
({\cal L}(\eth))_z\equiv (\Lam \ker(\eth^+_z))^*\otimes 
(\Lam \coker(\eth^-_z)) 
\] 
for each fixed $z\in B$ (as expected). 
If the family $\eth$ 
has index zero we have $\dim(H^+_{[0,\lambda)})= 
\dim(H^-_{[0,\lambda)})$ for each $\lambda>0$ and it makes 
sense to speak about $\det((\eth^+)_{[0,\lambda)}$ 
as a section  of $\det(H_{[0,\lambda)})$. 
These sections patch togeher (simply because 
$ 
\det(\eth^+_{[0,\mu)})=\det(\eth^+_{[0,\lambda)})\otimes 
\det(\eth^+_{[\lambda,\mu)}) 
$
and we obtain a smooth section $\det(\eth^+)\in\cin(B,{\cal L})$. 
This is the analogue of the section $\det(T)$ considered 
at the beginning of this section in the finite dimensional case.  
  
The definition of the determinant bundle involves only the {\it small} 
eigenvalues of the operators $\Lap^\pm_z$. This is not the case for 
the natural metric and metric-compatible connection, introduced by Quillen 
and Bismut and Freed respectively, which involve instead the 
full spectrum of $\Lap^\pm_z$. 
 
To define the  Quillen metric $\|\cdot\|_Q$ first observe that each 
$H_{[0,\lambda)}$, and thus each $\det(H_{[0,\lambda)})$, inherits a natural 
metric coming from the $L^2$-metric of $\cin(M^z,\Spi_z)$. The problem with 
this $L^2$-metric, which we denote by $|\cdot|_{\lambda}$, is that it is 
{\it not} well defined: there is discrepancy between $|\cdot|_{\lambda}$ 
and $|\cdot|_{\mu}$ equal to the product of the eigenvalues of $\Lap^+_z$ in 
the interval $(\lambda,\mu)$. Denote by $\zeta(s,\Lap^+_z,\lambda)$ 
the zeta function for the operator $P^+_{(\lambda,\infty)}(z) \Lap^+_z$. 
This gives a $\cin$ function $U_\lambda\ni z\rightarrow \zeta ' 
(0,\Lap^+_z,\lambda)$ and it is possible to see that the metrics 
\begin{equation} 
\|\cdot\|_{Q}=e^{-\zeta'(0,\Lap^+,\lambda)/2}|\cdot|_{\lambda} 
\label{eq:2.5.2} 
\end{equation} 
patch together to define a {\it global} metric on ${\cal L}(\eth)$. 
This is the Quillen metric; it involves the heat kernel  of $\Lap^+_z$ for 
all times. It gives another use for the determinant of a Laplacian as 
defined in \S 2.3; in fact for the section $\det(\eth^+) \in\cin(B,{\cal L})$ 
(which vanishes precisely when the operator $\eth_z$ is not invertible),
\[ 
\| \det(\eth^+) \|^2_Q = \det(\eth^-\eth^+). 
\] 
 
The Bismut-Freed connection is somewhat more complicated to describe and 
we shall not enter into the details here. Just like the Quillen metric, 
it is defined on each $\det(H_{[0,\lambda)})$ and then shown to be 
independent of choices. On each $U_{[0,\lambda)}$ the Bismut-Freed 
connection, henceforth denoted by $\nabla^{{\cal L}}$, is the sum of two pieces 
\begin{equation} 
\nabla^{{\cal L}}|_{U_\lambda}=\nabla^\lambda+\beta^+(\lambda). 
\label{eq:2.5.3} 
\end{equation} 
The first summand  $\nabla^\lambda$ is a connection which comes ultimately 
from the metric but is not globally defined; the second piece is a 1-form 
$\beta^+(\lambda)$ which  is given by a $t$-integral over $\RR^+$ 
involving $\eth^\pm$ and the heat-kernel $\exp(-t\Lap^\pm)$. This term 
should be thought of as a sort of eta invariant needed to make the
various definitions $\nabla^\lambda$ coherent. Bismut and Freed prove that 
this connection is compatible with the Quillen metric. Notice that
once again, even from this vague description, it is clear that we 
require the heat-kernel for all times. 
 
Given a family of Dirac operators $\eth=(\eth_z)_{z\in B}$ 
as above we now have a determinant line bundle ${\cal L}(\eth)$, 
with a natural metric $\|\cdot\|_Q$ and metric compatible 
connection $\nabla^{\cal L}$. One of the main contributions of Bismut-Freed is 
the explicit computatation of the curvature and holonomy of $\nabla^{\cal L}$; 
in other words they give geometric formulae for the local and global anomaly. 
We refer to their papers for a statement of the precise results. 
 
Our main concern here is of a different nature. Suppose as in \S 2.1 
that the fibration $M\rightarrow B$ defining the Dirac family is the 
union along a fibering hypersurface $H$ of two fibrations with 
boundary $M=M_+\cup_H M_-$. We have now four Dirac families, $\eth_M$, 
$\eth_{M_\pm}$ and $\eth_H$. If we fix a spectral section $P$ for 
$\eth_H$ then we obtain two families of Fredholm  operators,  
as in \S 2.1 and thus two determinant bundles ${\cal L}(\eth_{M_+},P)$, 
${\cal L}(\eth_{M_-},\Id-P)$. The questions we shall address later in 
the paper are
 
\begin{itemize}
\item Q1. Is there a natural isomorphism ${\cal L}(\eth)\longrightarrow 
{\cal L}(\eth_{M_+},P)\otimes{\cal L}(\eth_{M_-},\Id-P)$? 
\item Q2. Is it possible to define Quillen metrics and Bismut-Freed 
connections on these two line bundles, ${\cal L}(\eth_{M_+},P)$, 
${\cal L}(\eth_{M_-}, \Id-P)$, and prove surgery formul\ae\  for the 
corresponding curvature and holonomy? 
\end{itemize}

\section{A closer look at the methods of Bunke and Vishik}
Although the remainder of this paper is devoted mainly to a more 
detailed discussion of the surgery calculus of MM/HMM along with
a few of its applications, we wish to describe the other two 
principal methods, those of Bunke and Vishik, in at least 
a bit more detail than we have up until now.

\subsection{Bunke's unitary equivalence}
The method developed by Bunke \cite{Bu} was directed specifically 
at finding a gluing formula for the eta invariant. Continuing 
our discussion from \S 1.4, Bunke considers a compact manifold $X$ split 
along a hypersurface $H$ as usual, and with a metric $g$ containing 
an exactly cylindrical piece around $H$. On the manifolds with 
boundary, $X_{\pm}$, the Dirac operators are endowed with
augmented APS boundary conditions associated to a choice of
Lagrangian subspaces $\Lambda_{\pm} \subset \mbox{\rm ker\,}(\eth_H)$. 
The issue is to find a good expression for the defect 
\[
\delta = \delta(\Lambda_+,\Lambda_-) \equiv \eta(\eth_X) - 
\eta(\eth_{X_+},\Lambda_+) - 
\eta(\eth_{X_-},\Lambda_-).
\]

In the easier case, when $\eth_H$ is invertible, Bunke shows
that the reduction mod ${\ZZ}$ of the defect vanishes. He goes
on to obtain a formula for the (no longer reduced) defect 
in the general case:
\begin{equation}
\delta(\Lambda_+,\Lambda_-) = m(\Lambda_+,\Lambda_-) - 
2I(P_+,P_-) + \dim \mbox{\rm ker\,}(D_+) - \dim \mbox{\rm ker\,}(D_-).
\label{eq:bform}
\end{equation}
Here $m(\Lambda_+,\Lambda_-)$ is an invariant of the pair
of Lagrangian subspaces which is given in a few different ways.
The first is in terms of the averaged Maslov class, as explained
in \S 2.2, while the second, the form in which it was
originally found by Lesch and Wojciechowski \cite{LW}, is
as a sum of eigenvalues of sum matrix. This expression may
be written much more simply and neatly as in (\ref{eq:elexp}),
as discovered in \cite{asatet}. To explain the other terms
on the right, we must first explain his proof a bit more.

The core of the proof involves comparing two different Dirac operators. 
The first is the sum of Dirac operators $\eth_{X_\pm}$ on the disjoint 
sum of $X_\pm$, where these components are now assumed to have finite 
product cylindrical ends.  The other is for the sum of Dirac operators on 
the disjoint union of the manifold $X$, assumed to contain a long cylindrical 
piece around $H$, and the cylinder $C_R \equiv [-R,R] \times H$. On
each of these pieces, the Lagrangians $\Lambda_{\pm}$ are used
to augment the APS conditions at the appropriate boundary
components. These sums of operators are called $D_0$ and $D_+$, respectively. 
The operator $D_-$ is obtained from  $D_0$ by applying a unitary map, defined
using a partition of unity, which identifies the pieces of 
$X_+ \sqcup X_-$ with the equivalent pieces of $X \sqcup C_R$.
This is done only so that the operators $D_\pm$ live on
the same manifold. 

The other terms on the right in (\ref{eq:bform}) may now be
explained. The dimensions of the kernels of $D_{\pm}$ are the
obvious numbers. $P_{\pm}$ are the positive spectral projections
for the operators $D_{\pm}$, and after Bunke shows that their
difference $P_+ - P_-$ is compact, the relative index between
them, $I(P_+,P_-)$, is well-defined. 

In the nondegenerate case, only the first term on the right
in (\ref{eq:bform}) is necessarily trivial. However, if $\eth_X$ 
itself has only trivial nullspace (at least when the cylindrical 
piece is sufficiently long), then the other terms on the right
in (\ref{eq:bform}) also vanish. 

This formula is obtained by comparing the heat kernels of the 
operators $D_\pm$. This is accomplished by comparing a particular
regularization of the integrals required to define the
eta functions. These regularizations are
\[
R_{\pm}(s,t) \equiv \frac{1}{\Gamma(\frac{s+1}{2})}\int_0^\infty
r^{\frac{s-1}{2}}D_{\pm}e^{-(t+r)D_{\pm}^2}\,dr.
\]
The difference of eta invariants should arise as the limit
as $t \rightarrow 0$ of the difference of traces of the
$R_\pm$ at $s=0$. Unfortunately, these operators are not
continuous in the trace norm down to $t=0$, which makes
this procedure not entirely straightforward. Additional
terms are added on to ensure that the limit exists, and
these ultimately account for the various terms in the
expression (\ref{eq:bform}) for the defect. 

We note again that the dependence of the invariant
$m(\Lambda_+,\Lambda_-)$ on the Lagrangians $\Lambda_\pm$
was first considered by Lesch and Wojciechowski \cite{LW},
and they obtained one of the linear algebraic expressions
for this number, but not its identification with 
the averaged Maslov class. 

Amongst the advantages of this procedure are its relatively
elementary nature, and the reasonably explicit identification
of the integer part of the defect. 

\subsection{Vishik's variation of boundary conditions} 
The second approach, by Vishik, was developed for the study of 
determinants and analytic torsion. In \cite{V} Vishik studies
determinants of elliptic pseudodifferential operators in 
great generality and detail and gives, amongst other things, a 
new proof of the Cheeger-M\"uller theorem. Vishik's approach was 
recently adapted by Br\"uning and Lesch \cite{BL} to give another 
proof of the gluing formula for the eta invariant. Since this paper 
is somewhat more accessible than those of Vishik, and because we
have chosen to concentrate on the surgery formula for the
eta invariant specifically, we follow the discussion from this
latter paper instead.

Instead of considering a family of metrics degenerating 
(or lengthening) transversally to $H$, the perspective is 
now the `more classical' one, involving boundary conditions. 
The goal is to define a family of elliptic boundary problems 
$(\eth_X,\Pi_\theta)$ on the manifold $X_+ \sqcup X_-$. 
The boundary conditions are given by a family of orthogonal 
projections $\Pi_\theta$, $|\theta| < \pi/2$, acting on the 
direct sum $L^2(H;E) \oplus L^2(H;E)$. The two copies of this 
Hilbert space arise because $H$ needs to be thought of as the 
boundary of $X_+$ and of $X_-$ separately. To define $\Pi_\theta$
we use the following notation. For a (sufficiently smooth)
section $u$ on $X_+ \sqcup X_-$, denote its restriction to
$\del X_\pm$ by $u(\pm 0)$. Also, let $\Pi^{\pm}_{X_{\pm}}$
and $\Pi^0_{X_\pm}$ denote the spectral projections onto
the positive, negative and zero spectral subspaces of
$\eth_H$, where $H$ is considered alternately as the
boundary of $X_+$ and $X_-$, respectively. Then
\begin{eqnarray*}
\cos \theta \,\Pi^+_{X_+}u(+0) & = & \sin \theta \,\Pi^+_{X_-}u(-0) \\
\sin \theta \,\Pi^-_{X_+}u(+0) & = & \cos \theta \,\Pi^-_{X_-}u(-0) \\
\Pi^0_{X_+}u(+0) & = & \Pi^0_{X_-}
\end{eqnarray*}
defines the nullspace of the orthogonal projector on
$L^2(H;E)^2$ which we term $\Pi_\theta$. 

This family of boundary conditions interpolates between 
two extremes. One extreme, at $\theta = \pi/4$, is the 
`transmission condition': any section $u$ which solves 
$\eth u = 0$ on $X_+ \sqcup X_-$ and also $\Pi_0 u = 0$ 
must extend smoothly across $H$, and thus corresponds to an element 
of $\mbox{\rm ker\,}(\eth_X)$; the other, when $\theta = 0$,
corresponds to APS conditions independently on the
two pieces $X_\pm$. As the parameter $\theta$ changes,
the first projector `rotates' to the other.

Having defined these boundary conditions, one obtains
a family of self-adjoint problems, and for each operator
in this family one considers the eta invariant, which
we denote by $\eta(\eth_\theta)$. The main work in this
proof is showing first that the eta invariant for this
family of boundary problems is well-defined, i.e. that
the eta function is regular at zero, and then computing
the derivative of the eta invariant with respect to
$\theta$. After suitable normalizations, conjugating
by unitary transformations, Br\"uning and Lesch show
that this derivative vanishes. Unwinding the various
normalizations leads to the same gluing formula. 

While perhaps not quite as simple as Bunke's method, this 
method seems perhaps the simplest to adapt for the eventual 
study of gluing problems in more complicated situations. One
such situation arises when the hypersurface $H$ is the union 
of hypersurfaces with boundary $H_j$ intersecting at a common 
codimension two submanifold, $Y = \del H_j$ for all $j$. This 
generalization would be of particular importance if the signature 
formula of \cite{stmwc} is to be extended to manifolds with corners
of arbitrary codimension, cf. \S 5.3 below. 

\section{Pseudodifferential operators and the surgery problem} 

After the cursory treatments of the other methods in the previous
section, we now turn to a description of some of the details
of the surgery calculus of MM/HMM from \cite{asatei},
\cite{asatet}. We follow these papers, as well as \cite{tapsit} 
closely. In in the first two subsections below we continue and
extend the discussion of \S 1.3 on the $b$-calculus, and shall
refer to the notation there without further comment.
 
\subsection{Degenerating metrics} 
 
Let $X= X_+ \cup_H X_-$ as in the previous sections. Assume, 
just for the time being, that the Riemannian metric on $X$ is of 
product-type near $H$. Let ${\cal U}_H=(-T,T)\times H$ be a collar 
neighbourhood of $H$. Letting $T\rightarrow +\infty$ we obtain two 
manifolds, $\widehat{X}_+, \widehat{X}_-$, with infinite cylindrical 
ends. In general a manifold 
\[ 
\widehat{Z}=Z \cup_{\del Z} (\del Z\times (-\infty,0]) 
\] 
with a cylindrical end and metric $dt^2+h_{\del Z}$ along the cylinder  
$((-\infty,0])_t \times \del Z)$ can be compactified as a manifold 
with boundary $\overline{Z}$ with an exact $b$-metric, simply by 
making the change of variable $\log x=t$. In a neighbourhood
of the boundary $\{x=0\}$, corresponding to $t=-\infty$, the 
metric takes the form
\[ 
\frac{dx^2}{x^2}+h_{\del Z}. 
\] 
 
We can now give a natural analytic realization of the stretching 
procedure in this approach to the surgery problem. Let $h$ be an 
arbitrary  Riemannian metric on $X$ and let $x\in\cin(M)$ be a 
signed defining function for $H$; thus $H=\{x=0\}$, $dx\not=0$ on $H$,
but $x<0$ on $X_-$ and $x>0$ on $X_+$. Consider the one-parameter 
family of metrics 
\begin{equation} 
g_\e=\frac{dx^2}{x^2+\e^2}+h. 
\label{eq:4.1.1} 
\end{equation} 
For each  $\e>0$ the metric $g_\e$ is  non-degenerate on $X$. 
As $\e\rightarrow 0$ the metric $g_\e$ develops a neck across $H$,
the length of which is $L(\e) = (2\sinh^{-1}(\frac{1}{\e})+O(1))$.
The limiting metric, at $\e = 0$, is
\begin{equation} 
g_0=\frac{dx^2}{x^2}+h;
\label{eq:4.1.2} 
\end{equation} 
this is an {\it exact} $b$-metric on $X_+ \sqcup X_-$. 
We shall denote these two exact $b$-manifolds by $\overline{X}_+, 
\overline{X}_-$ and their disjoint union by $\overline{X}=\overline{X}_+ 
\sqcup \overline{X}_-$. In other words, $g_0$ endows the interior of 
$X_\pm$ with the structure of a Riemannian manifold with 
asymptotically (no longer necessarily product) cylindrical ends. 
 
In summary, the family of metrics (\ref{eq:4.1.1}) models the 
degeneration of $X$, through the stretching of a tubular 
neighbourhhod of $H$, from a closed compact manifold  to a 
manifold $\overline{X}$ which is the disjoint union of two 
manifolds with asymptotically cylindrical ends. 
 
As already explained the invariants we are interested in 
are each associated (one way or another) to the heat kernel 
of the Dirac Laplacian $\eth^2$. Since the heat operator is defined 
through the resolvent $(\eth^2-\lambda)^{-1}$ we conclude that  
the solution of the surgery problem for these invariants 
rests ultimately on a deep understanding of the uniform behaviour, 
down to $\e=0$, of the resolvent associated to $\eth_{X,g_\e}^2$, 
the Dirac Laplacian defined by metric (\ref{eq:4.1.1}). 
Amongst the other consequences of this analysis will
be a complete picture of the degeneration of the spectrum. 
 
The problem breaks into two intimately related problems. The first 
is to describe the {\it limit picture}, i.e. the geometry and 
analysis corresponding to $\e=0$. This is nothing more than the
setting of the $b$-calculus which we have already introduced. 
The second is to find a geometric setting (i.e. an appropriately
blown-up space) on which the relevant Schwartz kernels
behave uniformly with respect to $\e$. 
 
\subsection{The limit picture: more on the $b$-calculus} 
 
We have already discussed the $b$-calculus of pseudodifferential
operators on manifolds with boundary in \S 1.3. We shall now 
continue this discussion with two goals. On the one hand we shall
describe more precisely the analytic properties of the operators 
appearing in the limit of the surgery problem, while on the other 
hand parametrix construction for inverses of $b$-elliptic operators 
illustrates most of the main ideas in the more involved parametrix
construction in the surgery calculus. 
 
Getting down to business at last, consider the definition of the 
pseudodifferential calculus as given in \S 1.3. To motivate the 
introduction of the blown-up space there, we consider a simple example
drawn from \cite{tapsit} This is the $b$-differential operator 
\[ 
(1-x)x\frac{d}{dx}+c\quad c\in\RR 
\] 
on the manifold with boundary $Z=[0,1].$ This operator is invertible 
acting on $x^\alpha (1-x)^\beta H^1_b([0,1])$ if $\alpha<-c, \beta>c$; 
the Schwartz kernel of its inverse is given quite explicitly by the 
distribution 
\begin{equation} 
K_c(x,x')=-\frac{x^c}{(1-x)^c}\frac{(1-x')^c}{(x')^c}H(x'-x), 
\label{eq:4.2.1} 
\end{equation} 
where $H(\cdot)$ is the Heaviside function. This Schwartz kernel exhibits 
the usual diagonal singularities in the interior of $Z^2=Z\times Z$, 
and has additional singularities on the boundary hypersurfaces 
$\mbox{\rm lb}\equiv (\del Z\times Z)=\{x=0\}$, $\mbox{\rm rb}\equiv 
(Z\times \del Z)=\{x'=0\}$, and on the part of the corner which 
intersects the diagonal (viz. $(\del Z\times\del Z)\cap \Delta=\{(0,0) 
\cup (1,1)\}$). The corner carries the most complicated singularities, 
resulting from the interaction of those coming from the diagonal $\Delta$ 
and those coming from the two boundary hypersurfaces. The key observation 
is that the corner singularities can be simplified, or {\it resolved}, 
by introducing (generalized) polar coordinates. For example near 
$(0,0)$ consider the singular change of coordinates
\[  
r=x+x',\quad \tau=\frac{x-x'}{x+x'}. 
\] 
Then $2x=r(1+\tau), 2x'=r(1-\tau)$, and $K_c$ may be written as the 
product of a $\cin$ function and the distribution 
\begin{equation} 
\frac{(1+\tau)^c}{(1-\tau)^c}\times H(-\tau) 
\label{eq:4.2.3} 
\end{equation} 
which is singular on three {\it non-intersecting} hypersurfaces, 
$\tau=0$, $\tau=1$ and $\tau=-1$. These new variables are not smooth 
at the corner relative to the original ones, but on the $b$-stretched 
product $Z_b^2$ of \S 1.3 they are both smooth and independent.
Recall again that (assuming $\del Z$ is connected) $Z^2_b$ is the 
{\it blow-up} of $Z^2$ along $\del Z\times \del Z$, denoted $[Z^2 ; 
\del Z\times \del Z]$. If, as in the example above, $\del Z$ is not 
connected, then $Z^2_b=[Z^2;(\del Z\times \del Z)\cap\Delta]$. 
As a set $[Z^2 ; \del Z\times \del Z]$ is obtained by replacing the 
corner $\del Z\times \del Z$ with  its (inward pointing) spherical 
normal bundle $S_+ N(\del Z\times \del Z)$: 
\[ 
[Z^2 ; \del Z\times \del Z]=Z^2\backslash (\del Z\times \del Z)\sqcup 
S_+ N(\del Z\times \del Z). 
\] 
The $b$-stretched product comes with  a natural surjective 
blow-down map  
\[ 
\beta^2_b: [Z^2 ; \del Z\times \del Z]\longrightarrow Z^2 
\] 
and is given the minimal $\cin$ structure for which the lift 
of $\cin(Z^2)$ and of the polar coordinates around the corner  
are smooth. 
There are four important submanifolds in $Z^2_b$; the lifted diagonal 
$\Delta_b$, the left and right boundary faces, 
$\mbox{\rm lb},\,\mbox{\rm rb}$, obtained 
by lifting the corresponding boundary hypersurfaces in $Z^2$ 
and finally the new boundary hypersurface created by the blow-up: 
$\mbox{\rm ff}(Z^2_b)=S_+ N(\del Z\times \del Z)$. 
This is called the front face and by its very definition 
has the structure of a fibre bundle, with fibres diffeomorphic 
to the interval $[-1,1]$. 

The blow-up $[M;N]$ of a manifold with corners $M$ along the
the submanifold $N$ may be defined in substantial generality,
assuming only that $N$ satisfies certain local triviality
conditions. It formalizes the introduction of polar coordinates 
around $N$. We shall encounter other, more complicated, 
examples below. 
 
The $b$-calculus is defined by specifying the singularities 
allowed in the Schwartz kernels of its elements. As the example 
illustrates, and this is really the main point, 
these singularities are best understood when they are resolved, 
i.e. lifted to $Z^2_b$. 
 
In order to make the definition of \S 1.3 more precise recall first 
that if $M$ is a manifold with boundary $\{x=0\}$ and if $E\subset 
\CC\times\NN^+$ is a set of indices, then the space 
${\cal A}^E_{\mathrm phg}(M)$ of polyhomogeneous 
conormal functions can be defined. 
It consists of functions which are smooth in the interior and have an 
asymptotic expansion of the type 
\begin{equation} 
\sum_{(z,k)\in E} a_{z,k}x^z(\log x)^k,\quad a_{z,k}\in\cin(\del M) 
\label{eq:4.2.4} 
\end{equation} 
near the boundary. 
The index set $E$ specifies which exponents are allowed in  
(\ref{eq:4.2.4}); the set ${\NN} \times \{0\}$, which we denote 
by $0$ for brevity when the context is clear, corresponds to 
functions smooth up to the boundary. A similar definition may be 
given when $M$ is a manifold with corners; one simply requires
expansions of this type at all boundary hypersurfaces and product-type
expansions at the corners. Here it is necessary to specify 
an index family ${\cal E}=(E_1,\dots,E_n)$, where $E_j$ is an index set 
for the boundary hypersurface $H_j, j=1,\dots,n$. 
 
Now consider the case $M=Z^2_b$ and fix an index family 
${\cal E} = (E_{\rb},0,E_{\lb})$, where the boundary hypersurfaces
are listed in the order left boundary, front face and right
boundary. In particular, the index set associated to the front face 
$\ff(Z^2_b)$ is the one associated with smooth functions.
 
The $b$-calculus $\Psi^{*,{\cal E}}_b(Z)$ is the sum of two pieces: 
\[ 
\Psi^{*,{\cal E}}_b(Z)=\Psi^*_b(Z)+\widetilde{\Psi}_b^{-\infty,{\cal E}}
(Z). 
\] 
The first one, $\Psi^*_b(Z)$, is the {\it small} calculus: elements in 
$\Psi^*_b(Z)$ have Schwartz kernels on $Z^2$ lifting to 
$Z^2_b$ so as to have the usual interior singularities along $\Delta_b$, 
vanish to infinite order at $\lb, \rb$ and to be $\cin$ up to the front 
face $\ff(Z^2_b)$. To say that a conormal singularity along $\Delta_b$ is
smooth up to $\ff(Z^2_b)$ means that it extends smoothly across this
face as a distribution conormal to the extended diagonal.
For example, a $b$-differential operator has Schwartz kernel
which is a smooth delta section along $\Delta_b$, and hence
$\mbox{\rm Diff}_b^*(Z)\subset\Psi^*_b(Z)$ as expected.
The second summand contains the boundary terms, which are 
smooth in the interior and polyhomogeneous, with index family 
${\cal E}$ at the boundary of $Z^2_b$: 
\[ 
\widetilde{\Psi}_b^{-\infty,{\cal E}}(Z)=
{\cal A}^{\cal E}_{\mathrm phg}(Z^2_b).
\] 
Strictly speaking, it is also necessary to add to these
two summands a third, containing `very residual' parts
of the calculus. Since these are not important for the
present discussion, we shall not discuss them further. 

The small $b$-calculus is an algebra. The full $b$-calculus itself 
is not, but only for the trivial reason that sometimes the
boundary terms are not integrable. When they are, it is
possible to give precise composition formul\ae: when 
$A\in\Psi^{*,{\cal E}}(Z)$ and $B\in\Psi^{*,{\cal F}}(Z)$, 
the resulting operator $A\circ B$ will be an element of 
$\Psi^{*,{\cal G}}(Z)$, where the new index family ${\cal G}$ 
may be determined explicitly from ${\cal E}$ and ${\cal F}$. These 
composition formul\ae\ are one cornerstone of the whole theory; 
in some sense, the main work in setting up one of these
degenerate calculi is in proving such formul\ae.
They may either be proved directly, as in \cite{tapsit}, which
becomes less feasible in more complicated geometric situations,
or else using general facts about pushforwards of polyhomogeneous 
conormal distributions with respect to so-called {\it $b$-fibrations}
on manifolds with corners, cf. \cite{cocdmc}. 
Closely related to these arguments are those used to establish
the precise mapping properties for these operators, in particular 
their boundedness on weighted Sobolev spaces. 
 
For an invertible elliptic $b$-pseudodifferential operator $A$,
the inverse is an element of $\Psi_b^{*,{\cal E}}(Z)$ for 
some particular choice of the index set ${\cal E}$. This 
important result is, of course, the {\it raison d'\^{e}tre} for
establishing the calculus. 
It is not immediately apparent why it should be necessary 
to enlarge the small calculus to include the polyhomogeneous 
boundary terms. We explain this issue now. Let $P\in\mbox{\rm Diff}^m_b(Z)$ 
be a $b$-elliptic operator. Using a suitably adapted symbol 
calculus, we may construct a parametrix $Q_\sigma\in\Psi^{-m}_b (Z)$ for $P$.
This has the property that $P\circ Q_\sigma=\mbox{\rm Id}-R_\sigma$ with 
$R_\sigma\in\Psi^{-\infty}_b (Z)$. At this point it might seem that
we are essentially done, but this is not the case because 
elements of $\Psi^{-\infty}_b (Z)$ are {\it not} compact on $L^2$. 
In fact the  element $R_\sigma\in\Psi^{-\infty}_b (Z)$ is compact on $L^2$ 
if and only if the Schwartz kernel of $R_\sigma$ restricted to the front face 
is equal to zero. In one direction this property should be clear.
In fact, if $(K_{R_\sigma})|_{\ff}=0$, then $K_{R_\sigma}$ vanishes when restricted to 
any boundary face of $\del Z^2_b$, and so its pushforward to
$Z^2$ also vanishes on the entire boundary of $Z^2$. Compactness
of operators with Schwartz kernels of this form, which are
also smooth in the interior, follows from the Arzel\`a-Ascoli
theorem. 

To remedy this situation, we look for a correction term $Q'$ with the 
property that 
\begin{equation} 
P\circ (Q_\sigma - Q')=\mbox{\rm Id}-(R_\sigma - R')
\quad\mbox{\rm with}\quad (K_{R_\sigma})|_{\ff}= 
(K_{R'})|_{\ff}. 
\label{eq:4.2.6} 
\end{equation} 
Thus the operator $Q'$ is intended to cancel the restriction 
to the front face of $R_\sigma$. 
 
The restriction of the Schwartz kernel to the front frace is defined 
for any element in the small calculus; it defines a natural
homomorphism  
\[ 
I: \Psi^*_b(Z)\longrightarrow \Psi^*_b(\overline{N^+(\del Z)}). 
\] 
with $\overline{N^+(\del Z)}\cong [-1,1]\times \del Z$ the compactified 
inward pointing normal bundle. This map is called either the normal 
or indicial homomorphism, and we use these two names interchangeably.
(These two model operators exist for any of the degenerate calculi,
but coincide only in the special case of the $b$-calculus.) 
It is defined by observing that
the front face in $Z^2_b$ is canonically identified with the front 
face in the stretched product of $\overline{N^+(\del Z)}$.
Thus the restriction of the kernel to $\ff(Z^2_b)$ may be
transferred to the other stretched product and using the
dilation structure of $N^+(\del Z)$ may be extended further
to be homogeneous in the interior. This gives a kernel on 
$(\overline{N^+(\del Z)})^2_b$, i.e. a $b$-pseudodifferential operator  
on $\overline{N^+(\del Z)}$. In the special case of a $b$-differential
operator $P$, locally given by 
\[ P = \sum_{j+|\alpha| \leq m} 
a_{j,\alpha}(x,y)(x\del_x)^j\del_y^\alpha 
\] 
then it is not hard to check that this procedure leads to
the indicial operator for $P$,
\[ 
I(P)= \sum_{j+|\alpha| \leq m} 
a_{j,\alpha}(0,y)(x\del_x)^j\del_y^\alpha.
\] 
The normal homomorphism may be thought of as a noncommutative 
secondary boundary symbol.
 
To solve (\ref{eq:4.2.6}) then we need to find an operator $Q'$, 
defined by a Schwartz kernel in $Z^2_b$, such that 
$I(P)\circ I(Q')=I(R_\sigma)$. 
Formally 
\begin{equation} 
I(Q')=I(P)^{-1}\circ I(R_\sigma) 
\label{eq:4.2.6b} 
\end{equation} 
fixes the Schwartz kernel of $Q'$ near the front face, and this
may then be extended to all of $Z^2_b$. Formula (\ref{eq:4.2.6b}) shows 
that in order to construct a parametrix  we need to {\it invert} the 
normal operator of $P\in\mbox{\rm Diff}^m_b(Z)$. It is because
the inverse $I(P)^{-1}$ always involves polyhomogeneous boundary
terms that we must always include the polyhomogeneous part of the 
general $b$-calculus. This may be seen already in the one-dimensional 
example above. 
 
The invertibility of $I(P)$ is considered relative to weighted 
Sobolev spaces, and as in \S 1.3, the basic result is that 
except for a discrete set of values of the weight parameter 
$I(P)$ can be inverted; the inverse has polyhomogeneous expansions 
at $\lb$ and $\rb$. Different weights give rise to different index 
sets in the expansion. In the special case $\eth_X\in\Diff^{1}_{b}$,
the omitted set of weights coincides exactly to the spectrum of the 
boundary operator $\eth_{\del X}$; the elements in the various 
index families, i.e. the exponents allowed in the polyhomogeneous 
expansions, are given explicitly in terms of $\mbox{\rm spec}_{L^2}
(\eth_{\del X})$. The same sort of result also holds for $\eth^2$. 
 
In summary, we have indicated how, for each `admissible' weight 
$\delta$, i.e. one for which $I(P)^{-1}$ exists, this construction 
gives the Schwartz kernel of a right parametrix $G_\delta$ for 
$P\in\Diff^{m}_{s} (Z)$ acting on $x^\delta H^m_b$; 
$G_\delta$ itself is an element of $\Psi_b^{-m,{\cal E}(\delta)}$,
where the index family ${\cal E}(\delta)$ can be explicitly described. 
The remainder term $R_\delta = G_\delta P - \Id$ is compact
on $x^\delta L^2$. A left parametrix with similar properties
is constructed similarly. 

This parametrix construction may be applied to show that
the actual resolvent $(\eth^2-\lambda)^{-1}\in\Psi_b^{-m,
{\cal E}_\lambda}$, for some explicitly given index family 
${\cal E}_\lambda$. 
 
\subsection{The surgery calculus} 
Having described the limit picture for the surgery problem 
at $\e=0$ in the family of metrics (\ref{eq:4.1.1}), we now turn 
to a description of the uniform behaviour of the family of
resolvents $(\eth^2_{X,g_\e} - \lambda)^{-1}$. 
 
The basic idea is to incorporate the parameter $\e$ into
the geometric description of the Schwartz kernels. Thus consider
the space $M=X\times [0,\e_0]$, with projection $\pi_\e: 
M\rightarrow [0,\e_0]$. The metric $g_\e$ lifts to this space,
and is nondegenerate along the fibres of $\pi_\e$. The vector field 
$\sqrt{x^2+\e^2}\del_x$ is of (essentially) unit length 
with respect to this metric, and thus appears in the definition of  
the Dirac operator $\eth_{X,g_\e}$ (henceforth denoted 
simply by $\eth_\e$). This vector field is not smooth on $M$ --
it has singularities along the submanifold $H\times\{0\}=\{x=0,\e=0\}$. 
As usual, we resolve these singularities by blowing up
this submanifold.

We thus define the {\it single surgery space} $M_s$ as the blow-up 
of $M$ along $H\times\{0\}$: 
\[ 
M_s = [M;H\times \{0\}] = (M\backslash (H\times\{0\})\sqcup S_+N(H\times\{0\}). 
\] 
The single surgery space is equipped with a blow-down map $\beta_s: 
M_s\rightarrow M$ and thus with a projection $\pi_{s,\e}: M_s\rightarrow 
[0,\e_0]$. 
The set on the right hand side of the formula above is given the minimal 
$\cin$ structure containing both the lift of $\cin(M)$ 
and also the polar coordinate functions $(r,\theta)$ (with $x=r\cos\theta, 
\e=r\sin\theta$). 
 
Besides the uninteresting boundary at $\e=\e_0$ the single 
surgery space has two boundary hypersurfaces: the {\it $b$-boundary} 
$\overline{X}$, corresponding to the original boundary at $\e=0$, 
$\overline{X}=$(closure of $\beta^{-1}(\{\e=0\}\backslash H)$)= 
$\overline{X_-}\sqcup\overline{X_+}$, 
and the new boundary hypersurface created by the blow-up, 
the {\it surgery boundary} $\overline{H}=S_+ N(H\times\{0\})\cong 
[-1,1]\times H$. 
By construction, the singular vector field $\sqrt{(x^2+\e^2)}\del_x$ 
lifts to be smooth on $M_s$; in fact the lift belongs to  ${\cal V}_b(M_s)$. 
The latter space is the span, over $\cin(M_s)$, of 
lifts of the vector fields on $M$ that are tangent to $H\times \{0\}$. 
Clearly the Dirac operator $\eth_\e$ is in the algebra of operators
generated by vector fields in  ${\cal V}_b(M_s)$. 
However, $\eth_\e$ does not differentiate in the direction of 
the (lift of the) vector field $\e\del_\e$, and so we restrict 
our attention to the somewhat smaller class of {\it surgery
vector fields} on $M_s$, 
\[ 
{\cal V}_s(M_s)=\{V\in {\cal V}_b(M_s) : (\pi_{s,\e})_* V=0\}. 
\] 
The lift of the family of metrics $g_\e$ to $M_s$ is smooth 
and non-degenerate on ${\cal V}_s(M_s)$; moreover its restriction to 
$\overline{X}$ is precisely the exact $b$-metric $g_0$. The lifted metric 
may also be restricted to $\overline{H}$, and gives another 
exact $b$-metric there. 
 
The surgery differential operators on $X$, $\Diff^{*}_{s} (X)$, 
are now defined as the differential operators generated over $\cin(M_s)$ 
by the vector fields ${\cal V}_s(M_s)$. The notation 
$\Diff^{*}_{s}(X)$, referring to $X$ instead of $M_s$, is meant to 
indicate that these operators should be regarded as acting 
on $X$ (or rather, the fibres of $\pi_\e$, and depending parametrically 
in a precise manner on $\e$. The Dirac operator $\eth_\e$  
is a  surgery differential operators of order one. In fact it is 
{\it surgery-elliptic}, in the sense that may be locally
expressed by an elliptic combination of basis of sections of ${\cal V}_s(M_s)$. 
Similarly $\eth_\e^2\in\Diff^2_s$ and it is elliptic as well. (Henceforth
we shall merely write elliptic rather than surgery-elliptic). 
 
The main goal now is to define the surgery calculus, a pseudodifferential 
calculus naturally containing the inverses (when they exist) of the elliptic 
surgery differential operators. Surgery pseudodifferential operators are 
defined in terms of their Schwartz kernel on $X^2\times [0,\e_0]$. 
These are distributions on $X^2\times [0,\e_0]$ with specific 
singularities along the submanifolds $\Delta\times [0,\e_0]$, $H\times 
H\times \{0\}$, $H\times X\times \{0\}$, and $X\times H\times \{0\}$. 
It is convenient to introduce the notation 
\[ 
H_R=X\times H \quad\quad H_L=H\times X. 
\]
The Schwartz kernels of surgery operators are best described
as being pushed forward from the {\it surgery double space} $M^2_s$,
which is obtained from $X \times X \times [0,\e_0]$ by
blowing up these various submanifolds. The order in which
we perform these blow-ups is important. First we blow up 
$H\times H\times \{0\}$, obtaining the space $[X^2\times [0,\e_0]; 
H^2\times \{0\}]$ with its blow-down 
map $\hat{\beta}^2$. 
Then we blow up in $[X^2\times [0,\e_0]; H^2\times \{0\}]$ 
the lifts by $\hat{\beta}^2$ of $H_R\times \{0\}$ and $H_L\times \{0\}$. 
This defines $M^2_s$, and we denote this two-step blow-up 
process more succinctly by  
\[ 
M^2_s=[X^2\times [0,\e_0]; H^2\times \{0\}; H_R\times \{0\} \sqcup 
H_L\times \{0\}]. 
\] 
The total blow-down map is $\beta^2_s:M^2_s\rightarrow X^2\times [0,\e_0]$. 

As a general note about iterated blow-ups, the order in which
various submanifolds of a manifold with corners are blown up
is important and will in general affect the final space. There are 
various conditions on the submanifolds, however, which ensure that
the iterated blow-up may be performed in any order. 
 
The blow-ups in the definition of $M^2_s$ define three new boundary 
hypersurfaces. The blow-up of $H^2\times \{0\}$ produces the 
face $B_{\ds}$, and  the blow-up of $H_R\times\{0\}$ and $H_L\times \{0\}$ 
produces the hypersurfaces $B_{\rs}$ and $B_{\ls}$. We also have
the boundary hypersurface coming from the original boundary at $\e=0$
which is denoted by $B_{\db}$.  Finally the diagonal $\Delta \times [0,\e_0]$ 
lifts through $\beta^2_s$ to a submanifold $\Delta_s\subset M^2_s$. 
Notice that both $B_{\ds}$ and $B_{\db}$ have non-empty 
intersection with the lifted diagonal (the ``d'' in the subscript 
is meant to suggest this). 
 
The calculus of surgery pseudodifferential is the sum of two pieces: 
\[ 
\Psi^{*,{\cal E}}_s(X)=\Psi^*_s(Z)+\Psi^{-\infty,{\cal E}}_s(Z). 
\] 
The {\it small surgery calculus} $\Psi^*_s(X)$ consists of
operators with Schwartz kernels on $X^2\times [0,\e_0]$ which
are pushforwards from $M^2_s$ of distributions which 
exhibit the usual conormal singularities along the lifted 
diagonal $\Delta_s$ and which vanish to infinite order at the 
boundary hypersurfaces $B_{\rs}$ and $B_{\ls}$ (the ones not intersecting
the lifted diagonal). By construction, $\Diff^*_s(X)\subset \Psi^*_s(Z)$. 
 
The second piece of the calculus contains operators with nontrivial
boundary terms; their Schwartz kernels are smooth in the interior of 
$M^2_s$ but have polyhomogeneous conormal expansions of the type  
(\ref{eq:4.2.4}) at the various boundary faces. 
As in the discussion of the $b$-calculus, the exponents allowed 
in these expansions are given by an index family 
\[ 
{\cal E}=\{E_{\ds}, E_{\ls}, E_{\rs}, E_{\db}\} 
\] 
The boundary faces $B_{\ls}$ and $B_{\rs}$ are given index 
sets with strictly positive real part, which ensures that the
corresponding kernels vanish at these faces; the boundary faces 
$B_{\ds}$ and $B_{\db}$, the ones meeting $\Delta_s$, are given 
index sets with non-negative real part and with the first term in 
the expansion equal to $(0,0)$, which ensures that the
kernels can be restricted to these faces. It is possible to
consider index sets depending on a complex parameter, which 
will be the case for the resolvent $(\eth_\e^2-\lambda)^{-1}$, and
it also possible to discuss holomorphy in this context.

Again as with the $b$-calculus, it is most important to establish
how surgery pseudodifferential operators behave under
composition. There are composition formul\ae\ of the type
\[ 
\Psi^{m,{\cal E}}_s (X)\circ \Psi^{m',{\cal E}'}_s (X)\subset 
\Psi^{m+m',{\cal E}''}_s (X),
\] 
with ${\cal E}''$ given explicitly by ${\cal E}, {\cal E}'$. 
These are proved using the general results on pushforwards of
polyhomogeneous distributions in \cite{cocdmc}.
 
\subsection{The surgery resolvent} 
 
Having now defined the surgery calculus, one would like to show that 
the resolvent  $(\eth_\e^2-\lambda)^{-1}$ lies in it for a suitable 
choice of the index family ${\cal E} = {\cal E}(\lambda)$. 
This is proved by constructing a good parametrix for the
resolvent in this calculus, which we now sketch.

First consider the case where $\lambda\in\Omega$, $\Omega\cap [0,+\infty)
=\emptyset$. We wish to construct an element $E(\lambda)$ in the
surgery calculus which is an inverse of $(\eth_\e^2-\lambda)$ modulo 
a ``small'' remainder: 
\begin{equation} 
(\eth_\e^2-\lambda)\circ E(\lambda)=\mbox{\rm Id}-R(\lambda). 
\label{eq:4.4.1} 
\end{equation} 
Provided that the remainder term is sufficiently residual, the
right hand side of (\ref{eq:4.4.1}) can be inverted using
Neumann series, and after some work we can conclude that 
the resolvent itself is a surgery pseudodifferential operator. 

Using the symbol calculus, a version of which exists for
the small surgery calculus, we obtain an initial parametrix
$E_\sigma(\lambda) \in \Psi^{-2}_s(X)$, with 
\[ 
(\eth_\e^2-\lambda)\circ E_\sigma(\lambda)=\mbox{\rm Id}-R_\sigma(\lambda) 
\quad R_\sigma(\lambda)\in\Psi^{-\infty}_s(X). 
\] 
The remainder term $R_\sigma(\lambda)$ is compact when $\e > 0$,
but not when $\e = 0$. The problem comes from its nonvanishing 
restriction to the two boundary hypersurfaces meeting the lifted diagonal, 
$B_{\db}$ and $B_{\ds}$.

Exactly as we did in the $b$-calculus, we must then find a correction term 
which cancels the first term in the Taylor series of $R(\lambda)$ at 
$B_{\db}$ and $B_{\ds}$. In order to implement this argument, 
we use two {\it normal homomorphisms}, given by restrictions of
Schwartz kernels to these two boundary faces. To be more specific,
there are two natural identifications 
\begin{equation} 
B_{\ds}=[\overline{H}^2;\del\overline{H}^2]\quad  
B_{\db}=[\overline{X}^2;\del\overline{X}^2].
\label{eq:4.4.2} 
\end{equation} 
Restriction to each of these hypersurfaces defines, in turn,  
two surjective homomorphisms 
\begin{equation} 
N_s: \Psi^*_s(X)\rightarrow \Psi^*_b(\overline{H})\quad\quad 
N_b: \Psi^*_s(X)\rightarrow \Psi^*_b(\overline{M}), 
\label{eq:4.4.3} 
\end{equation} 
the surgery normal and $b$-normal homomorphism, respectively. 
 
Notice that in (\ref{eq:4.4.2}) we are blowing up the entire corner,
not just that component of it which intersects the lifted diagonal. 
The resulting $b$-calculi in (\ref{eq:4.4.3}) are therefore 
slightly larger than the ones considered in \S 4.2; the differences
in their properties, however, are negligible. 
 
These normal homomorphisms are also natural with respect to the
geometry. For example, 
\[ 
N_s(\eth_\e)=\eth_{\overline{H}}\quad\quad  
N_b(\eth_\e)=\eth_{\overline{X}}\equiv \eth_0,
\] 
where $\eth_{\overline{H}}$ is defined in terms of the restriction 
of the lift of $g_\e$ to ${\overline{H}}$.

Returning to the construction of a good parametrix, we must 
modify $E_\sigma(\lambda)$ by an operator $E(\lambda)' \in
\Psi^{-\infty,{\cal E}}_s$ for some index family ${\cal E}$
and such that 
\[ 
(\eth^2_\e - \lambda)\circ (E_\sigma (\lambda)-E(\lambda)')= 
\mbox{\rm Id}-(R_\sigma(\lambda)-R'(\lambda))\quad\mbox{with} 
\] 
\[ 
K_{R_\sigma(\lambda)}|_{\ds}=K_{R'(\lambda)}|_{\ds}\quad \quad 
K_{R_\sigma(\lambda)}|_{\db}=K_{R'(\lambda)}|_{\db}. 
\] 
This is equivalent to solving
\begin{equation} 
N_s(\eth_\e^2-\lambda)\circ N_s(E(\lambda)') = N_s(R_\sigma(\lambda)) 
\quad\quad 
N_b(\eth_\e^2-\lambda)\circ N_b(E(\lambda)') = N_b(R_\sigma(\lambda)). 
\label{eq:4.4.4} 
\end{equation} 
In other words, once again we need to {\it invert} the two normal 
operators: $(\eth_{\overline{H}}^2-\lambda)$ and $(\eth_0^2 - \lambda)$. 
For $\lambda$ away from the spectrum of $\eth_{\overline{H}}^2$ and 
$\eth_0^2 $, as we are at present assuming, this is possible. 
The solutions of these problems in (\ref{eq:4.4.4}) always
have nontrivial asymptotic expansions at $\del M^2_s$.
 
This argument fixes the Schwartz kernel of $E(\lambda)'$ on 
$B_{\ds}$ and $B_{\db}$ respectively, and we then find
some extension $E(\lambda)'$ to the entire space $M^2_s$.
The remainder term after this correction term has been added
is now sufficiently residual that we can iterate it away
without difficulty. 

In conclusion, we reemphasize that the fundamental step in
proving that the resolvent $(\eth^2_\e-\lambda)$, for $\lambda\in\Omega$, 
is an element of the surgery calculus $\Psi^{-2,{\cal E}_\lambda}_s (X)$ 
is the {\it inversion} of the two normal homomorphisms 
$N_s$ and $N_b$. 
 
\subsection{Small eigenvalues in the nondegenerate case} 
In analyzing the large time behaviour of the heat kernel 
$\exp(-t\eth_\e^2)$, uniformly in $\e$, it is necessary to 
understand the structure of the resolvent $(\eth^2_\e-\lambda)$ 
for $\lambda$ near $0$. The simplest case to understand
is when we impose the assumption that
\begin{equation} 
\mbox{the Dirac operator $\eth_H$ is invertible.} 
\label{eq:4.5.1} 
\end{equation} 
This hypothesis is called {\it nondegeneracy}. 
As indicated in \S 4.2, under this assumption the operator 
$\eth^2_0$ induced on $\overline{X}$ by the limiting 
metric $g_0$ is Fredholm on the ordinary (unweighted) $L^2$ space. 
In particular, $\mbox{\rm spec\,}(\eth^2_0)$ is discrete near $0$. 
This can be sharpened: if $\sigma_0^2$ is the smallest 
eigenvalue of the boundary operator $\eth_H$ associated to 
$\eth_0$, then the spectrum of $\eth_0$ is discrete in the 
interval $[0,\sigma_0^2)$. It is continuous, possibly with 
embedded discrete spectrum, in $[\sigma^2_0,+\infty)$. 
The full spectral and scattering theory of such operators
is described, using the $b$-calculus, in \cite{tapsit}.
A similar, but easier, analysis shows that assuming (\ref{eq:4.5.1}),
the surgery normal operator $N_s(\eth_\e)$ has only continuous 
spectrum contained in $[\sigma^2_0,+\infty)$. 
 
Choose $\delta$ so that $\mbox{\rm spec\,}(\eth_0^2)\cap (-\delta,\delta)
=\emptyset$. For $\lambda$ in a $\delta$-neighbourhood of $0$,
the resolvent of $\eth_0$ can be written as 
\[ 
(\eth_0^2-\lambda)^{-1}=\mbox{\rm Res}_0(\lambda)+\frac{1}{\lambda} 
\Pi_0,
\] 
where $\mbox{\rm Res}_0(\lambda)$ a parametrix depending holomorphically 
on $\lambda$ and with a finite rank error term, and 
$\Pi_0$ is the orthogonal projection onto the null space  
of $\eth_H^2$. We let $N=\mbox{\rm dim\,}\mbox{\rm null}(\eth_0)$. 
 
We can modify the construction of the resolvent in the surgery
calculus to take into account this refined structure of the inverse 
of the $b$-normal operator. In fact, it is not hard to produce
a surgery pseudodifferential operator $G(\lambda)$, depending 
holomorphically on $\lambda$ near zero, such that 
\begin{equation} 
(\eth_\e^2-\lambda)\circ G(\lambda)=\mbox{\rm Id}-\Pi(\lambda) 
\label{eq:4.5.4} 
\end{equation} 
with $\Pi(\lambda)$ a surgery pseudodifferential operator 
depending holomorphically on $\lambda$ and of uniform 
finite rank $= N$, and with $N_b(\Pi(\lambda))=\Pi_0$. 
Projecting the operator $\Pi(\lambda)$ onto its range, it is 
clear that the invertibility of the right hand side of 
(\ref{eq:4.5.4}) is equivalent to the invertibility 
of an $N\times N$-matrix of the form $(\delta_{ij}-a_{ij}(\lambda))$. 
If $q(\e,\lambda)$ is the determinant of this matrix, then 
it is holomorphic in $\lambda$ for each fixed $\e$ and 
polyhomogeneous conormal in $\e$. Moreover $q(0,\lambda)=
\lambda^N$. For each fixed $\e\in [0,\e_0]$, the function 
$q(\e,\cdot)$ has precisely $N$ zeros, counting multiplicity; 
these are the {\it small eigenvalues} of $\eth_\e^2$.
The orthogonal projection $\Pi_\e$ onto the small eigenvalues 
of $\eth^2_\e$ is therefore of uniform finite rank and it 
follows from this construction that it too is a surgery pseudodifferential 
operator. The small eigenvalues themeselves have polyhomogeneous expansion 
in $\e$.

In summary, assuming the nondegeneracy condition (\ref{eq:4.5.1}), 
for $\lambda$ in a small neighbourhood of zero, the resolvent 
$(\eth_\e^2-\lambda)$ is a meromorphic family of surgery pseudodifferential 
operators, with poles at 
the small eigenvalues of $\eth^2_\e$. The orthogonal projection 
onto the small eigenvalues is a surgery pseudodifferential 
operator of uniformly finite rank 
$N=\mbox{\rm dim\,}(\mbox{\rm null}(\eth_0))$.
 
\subsection{The logarithmic surgery calculus} 
The nondegeneracy condition (\ref{eq:4.5.1}) is strong, 
and often not satisfied in applications. To proceed further
without this assumption requires substantially more
work, unfortunately. Although we will not be able to
describe this in anywhere near the amount of detail
we have been going into up until now, we wish to
indicate a few of the new features in this general case.

The main problem is already seen in the very simplest
example of surgery degeneration, namely the one-dimensional
example of the interval $X = I = [-1,1]_x$ with the family of
metrics $g_\e = dx^2/(x^2 + \e^2)$. (The boundaries at
$x = \pm 1$ are unimportant here, and we could well have
considered the surgery degeneration of a circle at the
risk of slightly more complicated notation.) The total
length of $X$ with respect to $g_\e$ is $2L_\e$, where 
$L_\e = \mbox{\rm arcsinh\,}(1/\e)$, and the (Dirichlet)
eigenfunctions of the Laplacian $\Delta_\e$ are of the form
$u_k(r,\e) \equiv \sin(\pi k r/L_\e)$, $k \in {\ZZ}$, where 
$r = \mbox{\rm arcsinh\,}(x/\e)$. Already we see the new length 
scale: each of these quantities is most naturally expressed not 
in terms of the parameter $\e$, but rather in terms of the inverse 
logarithm of $\e$, $\ilg \e = 1/\log (1/\e)$. 

Proceeding further with this example, we next examine
the lifts of the eigenfunctions above on the single
surgery space $M_s$. After a brief calculation, we see
that $u_k$ lifts to a function obviously smooth in
the interior of $M_s$, equal to $1$ on the surgery
front face (the lift of $\{0\} \times \{0\}$), and
equal to $(-1)^k$ on the adjacent boundaries at $\e = 0$. 
The lift is not polyhomogeneous on $M_s$, and does not 
behave uniformly in $\e$. In fact, the oscillations of these 
eigenfunctions somehow disappear into the corners, at the intersection
of the surgery front face and the other $b$-faces at $\e = 0$. 

These various issues must be dealt with simultaneously,
and again the idea is to resolve these new singular
phenomena geometrically by performing some new blow-ups.
We describe these only for the singular surgery space,
and shall now define the single {\it logarithmic} surgery
space $M_{\Ls}$. The double logarithmic surgery
space $M_{\Ls}^2$ is unfortunately much more complicated
than the (already none-too-simple) space $M^2_s$,
and we must refer the interested reader to \cite{asatet}
for its definition.

To deal with the new length scale we first define
the logarithmic blow-up of $M_s$. This is obtained
by simply replacing the boundary defining function
$\rho$ of each boundary hypersurface (at $\e = 0$)
by $\ilg \rho$. In effect, this defines a new (but
equivalent) $\cin$ structure on $M_s$. Smooth
functions in this new structure are those which
are smooth in the various `new' functions $\ilg \rho$
on the original space. Although this may
not appear to be a blow-up in the sense we have
been describing this concept, it may be recast in
this language, cf. \cite{asatet}. Next we blow
up the corners, i.e. the intersections of the
$b$-faces and the surgery face at $\e = 0$. The
resulting space is now called the single
logarithmic surgery space $M_{\Ls}$. It has four
boundary faces at $\e = 0$, instead of the 
two possessed by $M_s$. 

Rather than entering into any more details of this
construction, suffice it to say that the overall
strategy is much the same as before. A double
logarithmic surgery space $M_{\Ls}^2$ is defined,
and is equipped with blow-down maps to $M_{\Ls}$. 
The space $\Psi_{\Ls}^{*,{\cal E}}$ of logarithmic
surgery pseudodifferential operators is again defined
as containing operators, the Schwartz kernels of which
are pushed forward from $M^2_{\Ls}$, and these kernels
on the double logarithmic surgery space are conormal
at all boundary faces (note that now this implies
the existence of expansions in powers of $\ilg \rho$
at any face with boundary defining function $\rho$).
The main theorem, proved by an explicit parametrix
construction, is that in the degenerate case the 
resolvent $(\eth_\e^2 - \lambda)^{-1}$ is an element of this
surgery calculus, in a precise sense uniformly even
as $\lambda$ approaches zero. 

The one aspect of this that we shall discuss slightly more 
is the new normal operator that must be considered. This is 
the {\it reduced normal operators}, $\mbox{\rm RN}(\eth_\e)$. 
To define it, first recall that since $0$ is now an eigenvalue
of $\ker(\eth_H)$, the limit operator $\eth_0^2$ is not Fredholm 
 on $L^2$ (the weight $0$ is one of the omitted ones). 
However, it is Fredholm on $x^{\pm\delta}L^2$ for $\delta$
sufficiently small, and has a parametrix $G(\pm\delta)$ in the
$b$-calculus. The structure of this parametrix can be
used to prove that near $H$  
\[ 
\eth^2_0 v=0, v\in x^{-\delta} L^2_b(\overline{X}_+ \sqcup 
\overline{X}_-)\Longrightarrow v\sim v_1\log x + v_0 +v' 
,\;\eth^2_H(v_i)=0\;,v'\in L^2_. 
\] 
These asymptotic boundary values define two pairs of subspaces 
in $\ker(\eth^2_H)$, analogous to the scattering Lagrangians
considered earlier for the Dirac operator. These are
\[ 
\Lambda^N_\pm=\{v_1 ; \exists\, v\sim  
v_1(y)\log x + v_0(y) +v', \eth^2_0 v=0, v'\in L^2\} 
\] 
\[ 
\Lambda^D_\pm=\{v_0 ; \exists \,v\sim v_0(y)+v', v'\in L^2\}. 
\] 
(The subscripts here refer to $X_\pm$). The reduced normal operator 
$RN(\eth_\e)$ is the operator $D^2_s$ on $[-1,1]$ acting on 
$\ker(\eth^2_H)$-values funtions with the boundary conditions 
\[ 
u|_{s=-1}\in\Lambda^D_- \quad D_s u|_{s=-1}\in\Lambda^N_- 
\] 
\[ 
u|_{s=+1}\in\Lambda^D_+ \quad D_s u|_{s=+1}\in\Lambda^N_+. 
\] 

Just as in the simpler parametrix construction in the 
nondegenerate case, we need to invert the various normal
operators, which now includes this new one. The
inversion of the reduced normal operator can be
done quite explicitly, and we can also see at least
in very vague outline how the scattering Lagrangians
enter into the analysis. 
 
\section{Applications of the surgery calculus} 
In this final section of this survey, we present four
applications of the surgery calculus. The first is
purely analytic, and is the detailed description
of how eigenvalues of $\eth_{X,\e}$ accumulate
as $\e \rightarrow 0$. The second is a final discussion
of the proof of the gluing formula for the eta
invariant from this point of view. 
After that we
discuss an interesting application of this gluing
formula, which is the signature formula for manifolds
with corners of codimension two. We finally discuss
 the  analytic torsion and  
some aspects of 
 the proof of the gluing formula
for determinant bundles. 

\subsection{Accumulation of eigenvalues} 
One of the immediate consequences of the construction
of the resolvent for the family of operators
$\eth_X^2$ is a formula for the rate of accumulation
of its eigenvalues as $\e \rightarrow 0$. We shall 
only state the result here and say very little 
about its proof, which involves the full intricacies of the 
logarithmic surgery calculus.

Consider the eigenvalues $\lambda_j(\e)$ of $\eth_X^2$. We are 
particularly interested in the eigenvalues tending to zero
as $\e \rightarrow 0$. In order to study them, we rescale 
by setting $\lambda_j(\e) = (\ilg \e) z_j(\e)$, where `$\ilg$' stands 
for the inverse logarithm, i.e. $\ilg \e = 1/\log (1/\e)$ as in \S 4.6. 
(Of course other eigenvalues $\lambda_j(\e)$ tend to finite nonzero
limits or to infinity, and it may be possible to study them
by analogous methods, but this has not been carried out.)
Recall also the reduced normal operator $\mbox{\rm RN\,}(\eth_\e^2)$
introduced at the end of the last section. 
\begin{theorem} 
Assuming that the eigenvalues $\lambda_j(\e)$, and hence $z_j(\e)$,
are listed in increasing order, with multiplicity, and similarly 
for the eigenvalues $\mu_j$ of the reduced normal operator 
$\mbox{\rm RN\,}(\eth_\e^2)$, then as $\e \rightarrow 0$, either
$z_j(\e) \rightarrow \mu_j$, or else $z_j(\e) \rightarrow \infty$. 
The number of eigenvalues converging to zero is the 
same as the dimension of the nullspace of the limiting operator 
$\eth_{X,0}^2$ on $X_+ \sqcup X_-$, and for each $\mu_j$, there is 
exactly one family of eigenvalues $z_j(\e)$ converging to it. 
\end{theorem}

The eigenvalues $z_j(\e)$ converging to $0$ are somewhat special: it 
can be proved that they are rapidly decreasing in $\ilg \e$ -- in fact,
they vanish as a power of $\e$ -- and are therefore called the {\it very 
small eigenvalues}. 

This result shows that the bottom of the spectrum of $\eth_{X,0}^2$ is 
somehow `granular', at least inasmuch as it is obtained as a limit 
of eigenvalues accumulating at a very slow rate.

The only point of the proof we wish to mention is
that it involves considering the resolvent with
rescaled spectral parameter
\[
R(z,\e) = (\eth_{X,\e}^2 - (\ilg \e)^2 z^2)^{-1}.
\]
The uniformity of this object is considered as
$\e \rightarrow 0$. The rescaling of the spectral
parameter corresponds essentially to blowing up
the $\lambda$-spectral plane at $\lambda = 0$. 

\subsection{The surgery formula for the eta invariant} 
We next consider the eta invariant for the Dirac operator associated 
to the metric $g_\e$: 
\[ 
\eta(\eth_\e)=\frac{1}{\sqrt{\pi}}\int^\infty_0 
t^{-1/2}\mbox{\rm Tr}(\eth_\e e^{-t\eth^2_\e})dt. 
\] 
According to the third approach to surgery, our 
main concern is to describe as precisely as possible 
the behaviour of $\eta(\eth_\e)$ as $\e\rightarrow 0$. 
 
First we assume that  
\begin{equation} 
\ker(\eth_H)=\{0\}. 
\label{eq:5.2.0} 
\end{equation} 
The large time behaviour of $e^{-t\eth^2_\e}$ can be analyzed 
using the results of the \S 4.5 : since the spectrum 
of $\eth^2_\e$ remains discrete near $\lambda=0$, down 
to $\e=0$, it follows from the contour integral representation 
\begin{equation} 
e^{-t\eth^2_\e}=\frac{i}{2\pi}\int_{\gamma}e^{-t\lambda} 
(\eth^2_\e-\lambda)^{-1} d\lambda 
\label{eq:5.2.1} 
\end{equation} 
(with $\gamma$ enclosing the spectrum) that $e^{-t\eth^2_\e}$ is 
exponentially decreasing as $t\rightarrow \infty$ uniformly 
in $\e$, up to a uniformly finite rank operator. 
 
Next one needs to understand the heat kernel uniformly for finite 
times. Consider first the heat kernel 
associated to a Dirac Laplacian on a closed Riemannian 
manifold $X$. For each $t>0$, $e^{-t\eth^2_X}$ is a smoothing 
operator. Bearing in mind the initial condition it must satisfy
we see that as $t\rightarrow 0$ the heat kernel must develop some 
sort of singularity along the diagonal. For the Laplacian in ${\RR}^n$, 
for example, the heat kernel is explicitly given by 
\[ 
\frac{1}{(2\pi t)^{\frac{n}{2}}} e^{-|x-x'|^2/4t}. 
\] 
In other words, viewed as a distribution on $X\times X\times [0,\infty)_t$, 
the heat-kernel is singular along the submanifold $\Delta\times \{t=0\}$. 
It is possible to encode the full short-time asymptotics of the
heat kernel, i.e. the expansion of this singularity, using 
the language of blow-ups; this is explained in detail in \cite{tapsit}. 
Because of the different homogeneities of the space and time variables, 
we must use {\it parabolic blow-up} instead of the normal blow-ups
to which we have mostly been referring. (In fact, these parabolic
blow-ups are akin to the logarithmic blow-ups of \S 4.6.)
The parabolic blow-up of a submanifold $Y$, which is defined
relative to a subbundle $S\subset N^*Y$, appeared first in \cite{EMM}. 
The {\it heat space} is the parabolic blow-up of $X^2\times 
[0,\infty)$ along $\Delta\times\{0\}$ with respect to the subbundle
$S=\mbox{\rm Span}\{dt\}$, and is denoted $[X^2\times [0,\infty);
\Delta\times \{0\},S]$. It is defined as the disjoint union of  
$ 
(X^2\times [0,\infty))\backslash (\Delta\times\{0\}) 
$ 
and a `parabolic normal bundle' to $\Delta\times\{0\}$. We shall not 
give more details of its definition, but only state the basic
fact that the fundamental solution of the heat equation lifts to 
a polyhomogeneous conormal distribution on this space. 

Because this construction is essentially local in the space
variables, it is also possible to define a $b$-heat space as
well as a {\it surgery heat space} $M^2_{\hs}$. This latter space
is the parabolic blow up of $M^2_s \times [0,\infty)_t$ 
along $\Delta_s\times \{0\}$ relative to the analogous subbundle
$S$ spanned by $dt$ along the diagonal at $\{t=0\}$. The heat
kernel of $\eth^2_\e$ is polyhomogeneous conormal on this
space, which fully encodes its uniformity as $\e \rightarrow 0$. 
In order to remain within the category of {\it compact} manifolds 
with corners, we can even compactify the temporal variable at $t=\infty$,
thus obtaining the {\it compactified surgery heat space} 
$M^2_{\mathrm c-hs}$. The notation $[0,1]_t$ denotes the
compactified $t$-axis, although of course $t$ is not the 
natural linear variable on this finite interval. The lifted diagonal 
embeds into $M^2_{\hs}$, and hence in $M^2_{\mathrm c-hs}$, and
this induces
\[ 
i_{\Delta_s}:\Delta_s\times [0,1]_t\equiv M_s\times [0,1]_t  
\hookrightarrow M^2_{\mathrm c-hs}. 
\] 

Using this notation, we can now reexpress the integral 
which define the eta invariant in terms of pull-backs
and push-forwards. Thus
\[ 
\eta(\eth_\e)=(\pi_t)_* (\pi_s)_* \left[i_{\Delta_s}^* \mbox{\rm tr}
\left(\frac{t^{-1/2}}{\sqrt{\pi}}\mbox{\rm Tr}(\eth_\e 
e^{-t\eth^2_\e})\right)\right]. 
\] 
Here $\mbox{\rm tr}$ denotes the trace on each fibre of the endomorphism 
bundle $\mbox{\rm Hom}(\Spi,\Spi)$. The map $\pi_s\equiv
\pi_s\times\Id :M_s\times [0,1]_t\rightarrow [0,\e_0]\times [0,1]_t$ 
is the composition of the blow-down map $\beta_s: M_s\rightarrow 
X\times [0,\e_0]$ and the natural projection $X\times [0,\e]\rightarrow 
[0,\e_0]$; finally $\pi_t:[0,\e_0]\times [0,1]_t \rightarrow  
[0,1]_t$ is the obvious projection. 

This formula expresses the eta invariant of $\eth_\e$ as the 
push-forward of the polyhomogeneous conormal distribution 
\[
i_{\Delta_s}^* \mbox{\rm tr}\left(\frac{t^{-1/2}}{\sqrt{\pi}}\mbox{\rm Tr}
(\eth_\e e^{-t\eth^2_\e})\right). 
\] 
from the manifold with corners $M_s\times [0,1]_t$ to the interval 
$[0,\e_0]$. Polyhomogeneity of the `integrand' at the two temporal 
faces follows from the short-time and large-time behaviour of the  
surgery heat kenel. The short-time behaviour follows from 
Getzler rescaling, as in \cite{tapsit}; it actually implies 
smoothness up to the face $t=0$. The large time 
behaviour has been already analyzed, and it implies the rapid 
vanishing at the corresponding temporal face, up to a uniformly finite 
rank operator. 
 
The resulting distribution on $[0,\e_0]$ may be analyzed using 
the general results on push-forwards from \cite{cocdmc}. 
In this manner, the behaviour of the eta invariant of $\eth_\e$ 
as $\e \rightarrow 0$ can be simply read off {\it geometrically}. 
To state the theorem recall that the projection $\Pi_\e$ 
onto the small eigenvalues of $\eth_\e$ is a finite rank operator, 
uniformly in $\e$. Let $\tilde{\eta}(\e)$ be the signature of 
$\Pi_\e$. Then we have  
\begin{theorem} (\cite{asatei}) 
The eta invariant associated to the Dirac operator 
$\eth_\e$ subject to the condition (\ref{eq:5.2.0}) satisfies 
\[ 
\eta(\eth_\e)={}^b\eta(\eth_{\overline{X}_+}) +
{}^b\eta(\eth_{\overline{X}_-})+ \tilde{\eta}(\e) + r_1(\e)+r_2(\e)\log\e
\] 
as $\e \rightarrow 0$, where $r_1, r_2 \in\cin([0,\e_0])$ with 
$r_1(0)=r_2(0)=0$. 
\end{theorem} 

When the operator $\eth_\e$ is no longer nondegenerate, a similar
proof works. One must define a logarithmic heat surgery
space, and then the eta invariant, as a function of $\e$, 
may be obtained as the push-forward from this space of 
a polyhomogeneous distribution. 
Let $\Pi_\e$ denote the orthogonal projection onto the eigenspaces
corresponding to the very small eigenvalues and let
$\tilde{\eta}(\e)$ be the signature of $\Pi_\e$. 
The generalization of
the previous result is
\begin{theorem} (\cite{asatet})
The eta invariant associated to the Dirac operator $\eth_\e$,
no longer necessarily satisfying the nondegeneracy hypothesis
(\ref{eq:5.2.0}), satisfies
\[
\eta(\eth_\e) = {}^b\eta(\eth_{\overline{X}_+}) + {}^b\eta(\eth_{
\overline{X}_-})+\tilde{\eta}(\e) + \eta(\mbox{\rm RN\,}(\eth_\e)) 
+ r_1(\ilg \e) + \log (\ilg \e) r_2(\ilg \e)
\]
as $\e \rightarrow 0$. Here, as before, $r_1$ and $r_2$ are smooth
functions vanishing at $0$. 
\end{theorem}
It is possible to calculate the eta invariant for the reduced
normal operator in terms of finite dimensional data involving
the scattering Lagrangian subspaces associated to the Dirac
operators on $X_\pm$, as discussed in \S 2.2. 

\subsection{The signature theorem on manifolds with corners} 
Suppose that $X$ is a compact manifold with corners. As with
manifolds with boundary, there are many possible choices
for natural metrics to consider on $X$. Following our
usual choice in this matter, we shall assume that the interior 
of $X$ is endowed with an exact $b$-metric $g$. This may
be described as follows. Assume that the codimension one
boundary faces of $X$ are listed as $M_\alpha$, $\alpha = 1,
\ldots N$, and that each such boundary face has a defining
function $x_\alpha$. Then near $M_\alpha$, 
\[
g = \frac{dx_\alpha^2}{x_\alpha^2} + h_\alpha,
\]
where $h_\alpha$ is some smooth nonnegative symmetric 2-tensor
in a collar neighbourhood of $M_\alpha$ which restricts to a metric
on $M_\alpha$, and near each corner of codimension $k$, given as 
the intersection of boundary faces $M_{\alpha_1}, \ldots, 
M_{\alpha_k}$,
\[
g = \frac{dx_{\alpha_1}^2}{x_{\alpha_1}^2} + \ldots + 
\frac{dx_{\alpha_k}^2}{x_{\alpha_k}^2} + h_{\alpha_1 \ldots \alpha_k},
\]
where the final summand restricts to a metric on the corner. 

Suppose that $\dim X = 4\ell$, and let $\eth_X$ denote the 
signature operator on $X$. This is simply the deRham-Hodge 
operator $d+ d^*$, restricted to act between the $+1$ and $-1$ 
eigenspaces of the natural algebraic involution $\tau$ which
equals $i^{p(p-1) + 2\ell}*$ on $p$-forms. The question
we discuss here is whether it is possible to obtain a
signature formula for $X$, relative to a metric of this 
(or any other) type. In the case where $X$ has only a boundary,
i.e. has corners only up to codimension one, then this
is precisely the celebrated signature formula for 
manifolds with boundary of Atiyah, Patodi and Singer \cite{APS}.
We have already mentioned some extensions and generalizations of 
this result, in particular its recasting by Melrose \cite{tapsit} 
and the families index theorems of \cite{MP2}, \cite{MP3}. In 
each of these papers, the signature formula is regarded as an 
index formula, and so it would seem that the most natural problem 
to study is whether it is possible to obtain an index formula 
for general Dirac-type operators associated to exact $b$-metrics 
on manifolds with corners. 

At this stage we recall for the reader the fact that on
a compact closed manifold there are in some sense two
index theorems: one for Dirac-type operators and
the other for general elliptic (pseudodifferential) operators. 
These are regarded as equivalent, because the latter may be 
deduced from the former using $K$-theory. One may consider these 
two types of index theorems for manifolds with boundary
or corners as well, but the relationship between them is no longer 
so simple. In this context the index theorem for general elliptic 
operators was only very recently obtained, by Melrose and Nistor 
\cite{MN1}, \cite{MN2}, but this formula is stated and proved using
Hochschild homology, and the terms in it do not translate readily to
more familiar ones for specific geometric operators. Thus it
is still an open problem to find an index formula for Dirac-type 
operators on manifolds with corners. 

Some partial progress on this sort of index theorem was made 
by M\"uller \cite{Mue2} when $X$ has corners only up to codimension 
two, assuming also some rather strong nondegeneracy conditions. If 
$\eth_X$ is a Dirac-type operator in the interior of $X$, 
then because of the nature of the metric, there are induced
Dirac-type operators $\eth_\alpha$ on every codimension one boundary 
face $M_\alpha$, each of which is now a manifold with boundary endowed 
with an exact $b$-metric, and also operators $\eth_{\alpha \beta}$
on the corners $H_{\alpha \beta} = M_\alpha \cap M_\beta$,
whenever these intersections are nontrivial. Since these
corners are compact, these latter operators have discrete
spectrum, but the $\eth_\alpha$ and $\eth_X$ have continuous spectrum. 
The continuous spectrum for the $\eth_\alpha$ is fairly simple, since 
it is of locally finite multiplicity, with thresholds at points 
determined by the eigenvalues of the $\eth_{\alpha \beta}$. In 
particular, $0$ is in the essential spectrum of $\eth_\alpha$ if 
and only if some $\eth_{\alpha \beta}$ is not invertible.
The continuous spectrum for $\eth_X$ itself is much more complicated, 
since its multiplicity is no longer necessarily locally constant. In 
particular, if some $\eth_{\alpha \beta}$ is not invertible, then the 
spectrum of $\eth_X$ near zero is of this rather complicated type.
In particular, it is unclear whether the basic method to understand 
this spectrum near zero, by analytically continuing the resolvent of 
$\eth_X^2$ to some branched cover of the plane, still works. 
Unfortunately, some sort of information about the spectrum of
$\eth_X$ near zero is necessary to obtain a formula for
the index, and this does not seem any too accessible. In any 
event, M\"uller's analysis assumes that each of the corner
Dirac operators is invertible so that, while the
continuous spectrum of $\eth_X$ may reach zero, it is
of the simpler type there, of locally finite multiplicity.
M\"uller does prove without this assumption, though, that 
when $\eth_X$ is the signature operator, then its $L^2$-index 
is well-defined and still yields the topological signature
of the manifold $X$. Unfortunately, the nondegeneracy
conditions are never satisfied for the signature operator,
and so M\"uller cannot deduce the signature formula
this way. 

It turns out that it is possible to obtain the signature
formula for manifolds with corners of codimension two
following a somewhat different sort of argument. This
was accomplished by the first author, Melrose and Hassell 
\cite{stmwc}. The idea is to obtain this formula not via an 
index calculation on all of $X$, but instead as a limit of 
index formul\ae\ on a family of compact manifolds with smooth 
boundary $X_\e$ which fill out $X$ as $\e$ tends to zero. 

There are a few steps to this proof. In the first, an
appropriate family of smoothings $X_\e$ is defined.
Then the APS signature theorem is applied to each of
these manifolds with boundary. Denoting by $\eth_{\e}$
the restriction of $\eth_X$ to $X_\e$, we get
\[
\mbox{\rm sign\,}(X) = \mbox{\rm ind\,}(\eth_\e)
= \int_{X_\e} \omega - \frac{1}{2} \eta(\eth_{\del X_\e})
+ B.
\]
The first term on the right here is the integral of the
usual signature density, which is the ${\cal L}$-polynomial
in the Pontrjagin forms, the second term is the eta
invariant of the induced signature operator on the
boundary, and the final term is an integral over $\del X_\e$ 
of a local expression involving the second fundamental form.
This final term is necessitated by the fact that the 
metric $g_\e$ is no longer of product type near the boundary.
The remainder of the proof involves calculating the
limits of these various terms as $\e \rightarrow 0$.
The left hand side, the signature, is topological, so
obviously does not change with $\e$. Using the asymptotics
of the metric $g$, the integral of $\omega$ over all of $X$
is well-defined, and the first term on the right tends
to this. The final term on the right tends to zero
because of the specific construction of the smooth
surfaces $\del X_\e$. 

Thus it remains to calculate the limit of the eta invariant. 
It turns out, again by the choice of smoothing, that the induced 
metric on $\del X_\e$ is simply undergoing surgery degeneration. 
Thus we already have the tools to analyze the limit of the eta 
invariant of the induced Dirac operator $\eth_{\del X_\e}$.
The only difficulty now is essentially combinatorial. The 
gluing formula for the eta invariant assumes only a single 
disconnecting hypersurface $H$, decomposing the manifold into 
two pieces. Here the relevant manifold $\del X_\e$ has some system 
of nonintersecting hypersurfaces $H_{\alpha \beta}$, and the 
metric is degenerating across each one of them. There are
again several ways of expressing the defect term in the
formula for the limit of the eta invariant. The most
elegant of these is as follows. Associate to $X$
a one-dimensional directed graph ${\cal G}$ by discarding
the interior of $X$, replacing each codimension one 
boundary component $M_\alpha$ by a vertex $v_\alpha$
and each codimension two corner $M_\alpha \cap M_\beta$
by an edge $e_{\alpha \beta}$. These edges are directed
by choosing arbitrarily some ordering of the $M_\alpha$,
then identifying $e_{\alpha \beta}$ in an orientation
preserving manner with $[-1,1]$ if $\alpha < \beta$ in
this ordering. We consider the trivial vector bundle
$V$ over ${\cal G}$ with fibre the direct sum of all
the cohomologies of the corners, i.e. the direct sum
of all $\mbox{\rm ker\,}(\eth_{\alpha \beta})$. There is
a Dirac operator acting on sections of this trivial
bundle. The `boundary conditions' at the vertex $v_\alpha$ 
are given by the scattering Lagrangian 
$\Lambda^\alpha_{\mathrm sc}$ associated to $M_\alpha$. The domain
of the Dirac operator $\eth_{{\cal G}}$ is the space of
sections $\phi$ which restrict along each edge $e_{\alpha \beta}$
to an element $\phi_{\alpha \beta}$ of the corresponding
nullspace $\mbox{\rm ker\,}(\eth_{\alpha \beta})$,
and such that at the vertex $v_\alpha$, the sum of
all $\phi_{\alpha \beta}$ for edges $e_{\alpha \beta}$
contiguous to that edge sum to an element of 
$\Lambda^\alpha_{\mathrm sc}$. Notice that in the simple case
where there are only two vertices and one edge, this
reduces to the operator introduced at the end of \S 2.2.
In any case, the defect term may be expressed as the
eta invariant of this operator 
$(\eth_{\cal G},\Lambda_{\mathrm sc})$. 
The final signature formula then is
\begin{theorem} Let $X$ be a manifold of dimension $4\ell$ with 
corners of codimension two, and suppose that $g$ is an exact 
$b$-metric on the interior of $X$. Then with the preceding
notation and conventions, 
\[
\mbox{\rm sign\,}(X) = \int_X \omega - \frac{1}{2}
\sum_{\alpha} \eta(\eth_{M_\alpha}) - 
\frac{1}{2} \eta(\eth_{\cal G},\Lambda_{\mathrm sc}).
\]
\end{theorem}
The correction term in this formula may once again be
expressed in purely finite dimensional linear algebraic 
terms using the various scattering Lagrangians, cf. 
\cite{stmwc}. 

The one place where we have really used special features
of the signature operator here is when we were able
to rule out any extra integer terms when taking the
limit of the eta invariant. This is because the
rank of $\eth_{\del X_\e}$ is determined topologically,
hence is constant. In general, there might well be
some spectral flow. The only thing we would be
able to deduce by this method in general, then, is
the mod ${\ZZ}$ reduction of this formula. 

\subsection{The surgery formula for the analytic torsion}
Behaviour of the analytic torsion under surgery was already studied 
in the fundamental work of Cheeger. Many of the subsequent proofs of 
the Cheeger-M\"uller theorem also exploit some form of this method 
in a basic way. In this section we look at the surgery problem 
for the analytic torsion from the point of view of the surgery
calculus, as studied by Hassell \cite{H}.

Let $X=X_+ \cup_H X_-$ and $g_\e=h + dx^2/(x^2+\e^2)$ be a family
of metrics on $X$ undergoing surgery degeneration along $H$, as usual.
Consider the analytic torsion $T(X,g_\e)$ associated to the metric 
$g_\e$. We can also consider the metric independent definition
given in \S 2.3
\[
T(X,\{\mu\})=T(X,g_\e)\cdot\Lambda(g_\e,\{\mu\})
\]
with $\{\mu\} = \{\mu^{(i)}\}$ a basis of $H^*(X) = \oplus H^i(X)$.
More generally, if $E$ is a flat unitary bundle we can define $T(X,E,g_\e)$
and $T(X,E,\{\mu\})$ using the de Rham complex twisted by $E$. 
A suitable understanding of the behaviour of the two terms appearing in 
this definition of $T(X,\{\mu\})$ will lead to a surgery formula
for $T(X,\{\mu\})$. The analysis of the first factor $T(X,g_\e)$ is based 
directly on the uniform analysis of the heat kernel associated to $\Lap_\e$, 
as in the case of the eta invariant. The construction of this heat kernel
as a polyhomogeneous distribution on the logarithmic surgery heat space is 
again the main ingredient in this analysis. The second factor 
$\Lambda(g_\e,\{\mu\})$ can be understood using a Hodge-theoretic
reinterpretation of the Mayer-Vietoris sequence for $X=X_+ \cup_H X_-$.

Putting these two results together Hassell proves \cite{H} that for 
suitable choices of $\{\mu\}$, 
\[
T(X,\{\mu\})={}^b T(\overline{X}_+,g_0) + {}^b T(\overline{X}_-,g_0)
+\ha\sum_{q=0}^n \,q\,\log\det RN(\Lap_q).
\]
The correct choice of set of bases $\{\mu\}$ is determined
by properties of the very 
small eigenvalues. These are rather simple to understand here 
because, using the (Hodge-) Mayer-Vietoris sequence again, 
it can be seen that the multiplicity of $0\in\mbox{\rm spec\,} 
(\Lap_\e)$ is {\it constant} in $\e \geq 0$. If $\Pi_\e$ is the 
orthogonal projection onto $\ker(\Lap_\e)$, then the $\{\mu\}$ 
in this formula must be chosen in the image of $\Pi_\e$, which is
by definition simply the Hodge cohomology of $(X,g_\e)$.

It is also important to use the fact that
\[
\ha\sum_{q=0}^n \,q\,\log\det RN(\Lap_q)
\]
may be explicitly decribed in terms of the finite dimensional 
subspaces $\Lambda^N_\pm$, $\Lambda^D_\pm$ appearing in the 
definition of the boundary condition for the reduced normal operator, 
and thus ultimately from the cohomology of $H$. In fact, 
another cohomological computation shows that this finite dimensional 
geometric expression also appears in a `surgery formula' for the 
combinatorial Reidemeister torsion $\tau(M,\{\mu\})$. 
Using this, it is possible to state the surgery formula for 
the analytic torsion in a particularly elegant way:

\begin{theorem} (\cite{H}) If $X=X_+ \cup_H X_-$ is odd dimensional,
the difference $\log T - \log\tau$ obeys the surgery formula
\[
\log\frac{T(X,g_\e)}{\tau(X,g_\e)}
=\log\frac{{}^b T(\overline{X},g_0)}{{}^b \tau(\overline{X},g_0)}
+\ha\chi(H)\log 2
\]
with $\overline{X}=\overline{X}_+\sqcup \overline{X}_-$
and $\chi(H)$ equal to the Euler characteristic
of $H$.
\end{theorem}
This result can be applied to reprove the Cheeger-M\"uller theorem
on the equality of the analytic and Reidemeister torsion on any 
closed compact manifold. Also, using a doubling argument it is also
possible to prove an extension of the Cheeger-M\"uller theorem for manifolds 
with boundary.

\begin{theorem} (\cite{H})
For an odd-dimensional manifold with boundary with exact $b$-metric $g$,
\[ 
{}^b T(Z,g)=2^{-\chi(\del Z)/4}
\,\tau(Z,g).
\]
\end{theorem}
Similar results hold when we twist by a flat unitary bundle $E$.

\subsection{Determinant bundles and surgery} 
In this section we finally address the two questions raised at 
the end of \S 2.5, following the treatment given by the second 
author in \cite{P2}, \cite{P3}. Recall the geometrical data: 
we are given a fibration $\phi:M\rightarrow B$ of compact 
manifolds with fibres even dimensional and  
endowed with smoothly varying metrics and smoothly varying 
spin structures. We denote by $g_{M/B}$ this family of fibre
metrics and by $\eth_M$ the associated family of Dirac operators.   
These data define a determinant bundle ${\cal L}(\eth)$ 
with a Quillen metric $\|\cdot\|_Q$ and Bismut-Freed connection 
$\nabla^{{\cal L}}$. We assume that the fibration $M$ is the union 
along a fibering hypersurface $H$ of two fibration with boundary: 
$M=M_+ \cup_H M_-$. These data fix the families $\eth_{M_\pm}$ 
as well as the family $\eth_{H}$. First let us make the very strong 
assumption that $\ker (\eth_H)_z = \{0\}$ for each $z\in B$. 
In this particular case, assuming that the metrics are product-like near 
$H$, the two families of APS boundary value problems on the fibrations 
$M_\pm$ are well defined and vary smoothly with $z\in B$; since
they are Fredholm they define two smooth determinant bundles, 
${\cal L}(\eth_{M_+},\Pi^+_0)$ and ${\cal L}(\eth_{M_-},\Pi^-)$, with 
$\Pi^+_0(z)$ equal to the spectral projection for $(\eth_H)_z$. 
 
These two determinant bundles can also be defined using the 
$L^2$ condition on the associated fibration with cylindrical ends:  
$\overline{M}=\overline{M}_+\sqcup\overline{M}_-$. 
In other words they can be defined in terms of the associated 
$b$-Dirac families $\eth_{\overline{M}_+}$, $\eth_{\overline{M}_-}$.  
We shall also use the suggestive notation $\eth_{\overline{M}} = 
\eth_{\overline{M}_-} \sqcup \eth_{\overline{M}_+}$. 
We denote the associated determinant bundles by ${}^b{\cal L}(\eth_{
\overline{M}_+})$ ${}^b{\cal L}(\eth_{\overline{M}_-})$. Notice that 
each Laplacian $(\Lap_{\overline{M}})_z$ has discrete spectrum near 
zero; this means that  the description of the $b$-determinant 
bundle in terms of small eigenvalues, as given in \S 2.5, 
is still valid. 
 
The determinant bundles in the two pictures (APS vs $L^2_b$) 
are canonically isomorphic (see \S 1.2); however there are substantial
advantages to working with $b$-determinant bundles. 
Namely, the definition of Quillen metric and Bismut-Freed connection 
can be given directly on ${}^b{\cal L}$, provided that the trace functional 
appearing in the definition of the zeta function $\zeta(s,\Lap^+,\lambda)
=\mbox{\rm Tr}(\Pi_{(\lambda,\infty)}(\Lap^+)^{-s})$ is replaced by 
the $b$-Trace and similarly for the second term, $\beta^+(\lambda)$, 
appearing in the Bismut-Freed connection (see \S 2.5). 
Here $\lambda$ must be always chosen away from the discrete spectrum 
of the family of $b$-Laplacians $\Lap_{\overline{M}}$. 
The first term $\nabla^\lambda$ in the definition of the Bismut-Freed 
connection is defined directly in terms of the metric and thus
extends to $b$-metrics with no effort. In summary, by using the 
$b$-Trace functional, we obtain in a natural way the $b$-Quillen metric 
$\|\cdot\|_{Q,b}$ and, more importantly, the $b$-Bismut-Freed connection  
\[ 
{}^b\nabla^{{\cal L}}|_{U_\lambda}={}^b\nabla^\lambda+{}^b\beta^+(\lambda). 
\] 
This latter step is not at all obvious in the APS framework. 
We note that in proving the compatibility of ${}^b\nabla^{{\cal L}}$ 
with $\|\cdot\|_{Q,b}$, the commutator formula for the $b$-Trace
is used in a crucial way.
 
Returning to the surgery problem, this discussion clarifies the 
{\it limit picture} at least under the assumption that 
$\Ker (\eth_H)_z = 0$. Thus let $x\in\cin(M)$ be a defining function 
for $H$ and consider the family of vertical metrics 
\[ 
g_{M/B}(\e)=\frac{dx^2}{x^2+\e^2}+g_{M/B}. 
\] 
Let $\eth_{M}(\e)$ be the associated Dirac family on the closed fibration 
$(M\rightarrow B, g_{M/B}(\e))$. We denote by $\nabla^{{\cal L},\e}$ 
the associated Bismut-Freed connection. 
We denote by ${}^b\nabla^{{\cal L}}_+$ and ${}^b\nabla^{{\cal L}}_-$ 
the $b$-Bismut-Freed connections induced by the limit metric 
$g_{M/B}(0)$ on the fibrations $\overline{M}_+$, 
$\overline{M}_-$. 
 
The following theorem is proved in \cite{P2} under the assumption 
$\Ker (\eth_H)_z = 0$ for all $z\in B$.
\begin{theorem} 
There exists a natural explicit isomorphism of determinant bundles 
\[ 
S(\e): {\cal L}(\eth_M (\e))\longrightarrow {}^b{\cal L}(\eth_{\overline{M}_+}) 
\otimes {}^b{\cal L}(\eth_{\overline{M}_-}). 
\] 
For the curvature and the holonomy of the corresponding Bismut-Freed 
connection the following formul\ae\ hold: 
\[ 
\lim_{\e\rightarrow 0}\,\,(\nabla^{{\cal L},\e})^2 =  
({}^b\nabla^{{\cal L}}_+)^2 
+({}^b\nabla^{{\cal L}}_-)^2 
\] 
\[ 
\lim_{\e\rightarrow 0}\,\,\mbox{\rm hol}_\gamma (\nabla^{{\cal L},\e})= 
\mbox{\rm hol}_\gamma ({}^b\nabla^{{\cal L}}_+)\cdot 
\mbox{\rm hol}_\gamma ({}^b\nabla^{{\cal L}}_-) 
\quad 
\forall\gamma\in\mbox{\rm Map}(S^1,B). 
\] 
\end{theorem} 
The proof is another application of the surgery calculus; the 
explicit isomorphism is induced by the projection $\Pi_\e$ onto 
the {\it small} eigenvalues of $\Lap^\pm$. The behaviour of the 
curvature and the holonomy of the Bismut-Freed connection is obtained 
by working directly with a push-forward of the latter object. 
The first ``metric'' part of $\nabla^{{\cal L},\e}$ (see 
(\ref{eq:2.5.3})) converges almost by definition; the second part, 
i.e. the term $\beta^+_\e(\lambda)$ can be analyzed using the 
heat-surgery calculus. One can prove that for $\lambda$ small  
\[  
\lim_{\e\rightarrow 0} \, \beta_\e^+(\lambda) + \log\e\cdot 
d\zeta'(0,\eth^2_H,0)= 
{}^b\beta^+_{\overline{M}_+}+{}^b\beta^+_{\overline{M}_-} 
\] 
where the convergence is to be taken as $C^k$ convergence 
of 1-forms on the set $U_\lambda$. The two formul\ae\ in the 
theorem then follow readily.
 
This theorem successfully solves the surgery problem on determinant 
bundles under the nondegeneracy condition on $(\eth_H)_z,\,z\in B$. 
We now drop this assumption and consider the general case 
(see \cite{P3}). Since $\eth_H$ arises as a boundary family we can 
certainly fix a spectral section $P$ for $\eth_H$ obtaining as in 
\S 2.1 and \S 2.5 the two Fredholm families $(\eth^+_{M_+},P)$, 
$(\eth^+_{M_-},\mbox{\rm Id}-P)$ and thus  
the two determinant bundles ${\cal L}(\eth_{M_+},P)$, 
${\cal L}(\eth_{M_-},\mbox{\rm Id}-P)$. 
These should be thought of as APS-determinant bundles. 
However, in order to apply the surgery calculus and the $b$-calculus 
we need to consider the corresponding $b$-determinant bundles, as we did 
in the invertible case.  This is indeed possible provided 
we allow the use of pseudodifferential operators. 
 
We shall now briefly pause to explain this fundamantal point. 
Let $\Dir=(\Dir_z)_{z\in B}$ be any family of Dirac operators  
on manifolds with boundary. Let us denote by $\Dir_{\del}$ the 
boundary family and let $P$ be a spectral section for $\Dir_{\del}$. 
It is proved in \cite{MP1} that  there exists a smooth family of 
operators $A^P$, with $A^P(z)\in\Psi^{-\infty}_b$ with the following 
properties: 
\begin{itemize} 
\item The family $(\Dir+A^P)_{\del}$ and the family of indicial operators 
$I(\Dir+A^P)$ are both  invertible; 
according to the $b$-calculus the family $(\Dir+A^P)$ is  
then Fredholm on $L^2$. 
\item The two families of Fredholm operators, one defined by the 
generalized APS-boundary value problem $(\Dir,P)$, and the other 
fixed by the $b$-family $(\Dir+A^P)$, are {\it homotopic}. 
\item The family $(A^P)_{\del}$ is finite rank and self-adjoint,
and $P_z$ is equal the projection onto the 
non-negative part of the spectrum of $((\Dir+A^P)_{\del})_z$ 
for each $z\in B$. 
\end{itemize} 
We refer to $A^P$ as a $P$-trivializing perturbation. 
These properties establish the following important priciple: 
the general APS family index theory defined by a spectral section $P$ 
can always be reduced to the {\it invertible} case but only by
passing to a {\it larger} class of operators. In any case, using
these properties, we now have a $b$-determinant bundle defined 
in terms of the Fredholm family $(\Dir+A^P)$, which we denote by 
${}^b {\cal L}(\Dir+A^P)$. Since we are again in the
invertible case, we can use the $b$-calculus and introduce 
a $b$-Quillen metric and a $b$-Bismut-Freed connection, 
essentially  as in the previous case.
 
Returning again to the surgery problem, this discussion
shows that there are $b$-determinant bundles, ${}^b {\cal L}(\eth_{
\overline{M}_+}+A^P_+)$ and ${}^b {\cal L}(\eth_{\overline{M}_-}+
A^{(\id-P)}_-)$ endowed with $b$-Quillen metrics and Bismut-Freed 
connections, ${}^b\nabla^P_+$, ${}^b\nabla^{(\id-P)}_-$. 
The surgery calculus can be used to show the existence 
of an element in the (fibre) surgery calculus $A(\e)\in\Psi^{-\infty}_s$ 
with the property that as $\e\rightarrow 0$ 
\[ 
\eth(\e)+A(\e)\longrightarrow (\eth_{\overline{M}_+}+A^{P}_+) 
\sqcup (\eth_{\overline{M}_-}+A^{(\id-P)}_-) 
\] 
(in the precise sense of \S 4.3). 
Since the family $\eth(\e)+A(\e)$, $\e>0$ is a perturbation by
a family of smoothing operators of the family $\eth(\e)$, it is 
certainly Fredholm. The associated determinant bundle ${\cal L}(\eth(\e)
+A(\e))$ can be endowed with a Quillen metric and Bismut-Freed connection. 
The arguments leading to the theorem above can now be extended 
(using the full force of the surgery and $b$-pseudodifferential 
calculi), resulting in the explicit isomorphism 
\[ 
S_P(\e): {\cal L}(\eth(\e)+A(\e))\longrightarrow  
{}^b{\cal L}(\eth_{\overline{M}_+}+A^{P}_+)\otimes 
{}^b{\cal L}(\eth_{\overline{M}_-}+A^{(\id-P)}_-) 
\] 
and the (asymptotic) additivity of the curvatures and multiplicativity 
of the holonomies. 
 
The final step is to show that these surgery formul\ae\ 
for the curvature and the holonomy are {\it independent} of 
the particular choice of perturbations $A^{P}_+ \sqcup A^{(\id-P)}_-$ 
and $A(\e)$. 
 
Consider first the closed case. Let $\eth$ be a Dirac family and  
$A^0$, $A^1$ two {\it smoothing} perturbations. We obtain two 
determinant bundles, ${\cal L}(\eth+A^0)$ and ${\cal L}(\eth+A^1)$, 
endowed with their hermitian structures. The space of smoothing 
perturbations is clearly simply connected; let $A(r), r\in [0,1],$ be 
a path of perturbations. Consider the family $\DDir$ on $B\times [0,1]$ 
given by $(\DDir)_{(z,r)}=\eth_z + (A(r))_z$. This is a family 
of Fredholm operators and we can consider the associated 
determinant bundle. The latter is endowed with a Bismut-Freed connection 
$\nabla^{\DDir}$; the local anomaly formula of Bismut-Freed can now 
be applied, and it shows explicitly that the curvature of $\nabla^{\DDir}$ 
is zero in the $dr$-direction. Since the space of smoothing perturbations
is simply connected, this shows that parallel transport defined by 
$\nabla^{\DDir}$ gives a {\it canonical} isomorhism 
$\tau: {\cal L}(\eth+A^0)\rightarrow{\cal L}(\eth+A^1)$  
which preserves curvature and holonomy. This property can be applied 
to the pair  of families $\eth(\e)$ and $\eth(\e)+A(\e)$ 
as well as the pair $\eth(\e)+A(\e)$ and $\eth(\e)+B(\e)$  
for a different choice of perturbation $B(\e)$. 
 
Consider now the boundary case and, as above, let $\Dir=(\Dir_z)_{z\in B}$ 
be a family of Dirac operator on manifolds with boundary. 
Denote by $\Dir_{\del}$ the boundary family and let $P$ 
be a spectral section for $\Dir_{\del}$. It is not difficult to see, 
using the $b$-calculus,  that the space of $P$-trivializing 
perturbations is simply connected. Let $A_0^P$ and $A_1^P$ two 
$P$-trivializing perturbations, $A^P(r)$ a path joining them and 
$\DDir_P$ the induced family on $B\times [0,1]$. 
The Bismut-Freed curvature formula is extended to manifolds with 
boundary in \cite{P3}. Applying the formula we discover that the 
$dr$-component of the curvature of the $b$-Bismut-Freed connection  
of the determinant bundle associated to $\DDir_P$ is {\it not} zero; 
however, and this is the key point, {\it it only depends on  
the boundary family $(\DDir_P)_{\del}$}. 
When this argument is applied to the fibration 
$\overline{M}=\overline{M}_+ \sqcup \overline{M}_-$ 
the two contributions cancel out because of the different 
orientation of the normals; thus the parallel 
transport defined by the Bismut-Freed connection 
produces, as in the closed case, a {\it canonical} isomorphism preserving 
curvature and holonomy. This means that the surgery results  
established for the $P$-trivializing perturbation $A^{P}_+ \sqcup 
A^{(\id-P)}_-$ on $\overline{M}= \overline{M}_+ \sqcup \overline{M}_-$ 
and the surgery perturbation $A(\e)$, only depend on the family $\eth_M(\e)$,  
the limit families $\eth_{\overline{M}_\pm}$ 
and on the choice of spectral section $P$ for the family of 
operators $\eth_H$ induced on the fibering hypersurface  
defining our decomposition $M=M_+ \cup_H M_-$. 
These results answer the questions raised at the end of \S 2.5 
in the framework of the $b$-calculus. It is still an open problem as to 
whether the APS-framework, and the other two approaches to surgery, 
can be used to give similar answers.

\end{document}